\documentclass[%
reprint,
preprintnumbers,
superscriptaddress,
amsmath,amssymb,
aps,
prl,
onecolumn
]{revtex4-2}
\usepackage{subcaption}    
\usepackage{booktabs}
\usepackage{overpic}
\usepackage{graphicx}
\usepackage{bm}
\usepackage[usenames,dvipsnames,table]{xcolor}
\usepackage[colorlinks=true,linkcolor=Blue,citecolor=Blue,urlcolor=Blue]{hyperref}
\usepackage[capitalise]{cleveref} 
\crefname{section}{Sec.}{Secs.}

\usepackage{mathtools}
\usepackage[justification=raggedright, singlelinecheck=false, font=small]{caption}

\DeclareMathOperator{\diag}{diag}

\newcommand{\eps}{\varepsilon}

\newcommand{\abs}[1]{\left\lvert #1 \right\rvert}

\newcommand{\ord}[2]{#1^{(#2)}}

\usepackage{tikz}
\usetikzlibrary{shapes.geometric, arrows.meta, positioning}

\begin{document}

    \preprint{RIKEN-iTHEMS-Report-26}

	\author{Riccardo Muolo}
    \email{riccardo.muolo@riken.jp}
	\affiliation{RIKEN Center for Interdisciplinary Theoretical and Mathematical Sciences (iTHEMS), Saitama, Japan}
	\affiliation{Department of Systems and Control Engineering, Institute of Science Tokyo, Tokyo, Japan}
	
	\author{Hiroya Nakao}
	\affiliation{Department of Systems and Control Engineering, Institute of Science Tokyo, Tokyo, Japan}
	\affiliation{Research Center for Autonomous Systems Materialogy, Institute of Science Tokyo, Yokohama, Japan}
	
	\author{Christian Bick}
    \affiliation{Department of Mathematics, Vrije Universiteit Amsterdam, Amsterdam, The Netherlands}
    \affiliation{Institute for Advanced Study, Technical University of Munich, Garching bei München, Germany}
    \affiliation{Department of Mathematics, University of Exeter, Exeter, United Kingdom}
    \affiliation{Mathematical Institute, University of Oxford, Oxford, United Kingdom}

    \title{Physical and emergent nonpairwise interactions in oscillator networks: from higher-order phase reduction to coupling design}
    
	\date{\today}
	
	\begin{abstract}
		 Phase reduction is a powerful technique to obtain phase models from high dimensional oscillatory systems. Starting from pairwise interactions, a first order approximation yields Kuramoto and Winfree phase models, while nonpairwise interactions emerge at the second order. Recently, phase models have been extended to account for nonpairwise interactions, notably, the higher-order Kuramoto model. Is there an intrinsic difference between physical and emergent nonpairwise interactions? And can we make use of the former to even out the effects of the latter? In this work, we exploit a recently developed parametrization method to compute the phase reduction and answer these questions. After a comparison of the new method with the classic Kuramoto-style phase reduction, we derive the network motifs of the emergent nonpairwise interactions and compare them with the physical once. Lastly, we exploit the motifs to adopt a coupling design approach. Our framework paves the way for further exploitations of physical nonpairwise interactions for applications in synchronization engineering.
	\end{abstract}
	
	\maketitle

\section*{Introduction}
\label{sec:introduction}

Synchronization is a fascinating behavior, consisting of the spontaneous emergence of coherent dynamics in populations of interacting self-sustained oscillators. 
This phenomenon has been observed since the classical experiments of Huygens in the 17th century~\cite{pikovsky2001synchronization}, and it appears in many natural and engineered systems, including neuronal networks, power grids, chemical oscillators, and mechanical systems~\cite{pikovsky2001synchronization,strogatz2000kuramoto,arenas2008synchronization}. 
A common theoretical approach is to represent each unit as an oscillator and to study how interactions between oscillators lead to collective dynamics. 
Understanding which types of interactions promote or hinder synchronization is of paramount importance in applications, as synchronization can be both beneficial (e.g., heart, power grids) or detrimental (e.g., mechanical systems). Moreover, the possibility of designing the interactions between oscillators to obtain a desired collective behavior is at the basis of synchronization engineering~\cite{kori2008synchronization,kiss2018synchronization}.

Recent works have highlighted the importance of nonpairwise interactions in nonlinear dynamics~\cite{bick2023higher,battiston2020networks,carletti2020dynamical,carletti2020random,bianconi2021higher,natphys,majhi2022dynamics,boccaletti2023structure,muolo2024turing,millan2025topology,battiston2026collective} and, in particular, how they affect synchronization (see, e.g.,~\cite{krawiecki2014chaotic,bick2016chaos,bick_explosive,gambuzza2021stability,gallo2022synchronization,della2023emergence,nurisso2024unified,muolo2025pinning,zhang2024deeper,skardal2025mixed,muolo2026higher,wang2026moderate,moriame2026efficiency}). 
These interactions involve three or more units simultaneously and cannot be decomposed into sums of pairwise contributions for a given choice of coordinates~\cite{Aguiar2020,neuhauser2020multibody,Muolo2023turing}. 
In the literature, \textit{nonpairwise} interactions are often called \textit{higher-order} interactions\footnote{They can also be called \textit{polyadic} interactions, \textit{many-body} interactions, \textit{group} interactions, etc.}. In some cases, this terminology can be ambiguous, as \textit{higher-order} can also denote terms in a perturbative expansion. Since in this work we consider phase reduction, i.e., a perturbative expansion with respect to the coupling strength, we will always use the term \textit{nonpairwise} to denote such interactions.

Where nonpairwise interactions arise depends on the level of description and on the model of the physical system. If the states in the physical system are given by~$\vec{x}_k,\vec{x}_j,\vec{x}_l\in \mathbb{R}^d$, then a nonpairwise interaction may correspond, for example, to a nonlinear coupling function~$\vec{G}(\vec{x}_k,\vec{x}_j,\vec{x}_l)=(x_{k1}x_{j1}x_{l1}, \dotsc, x_{kd}x_{jd}x_{ld})^\top$. Such terms appear naturally in systems with nonlinear responses or mediated interactions. If the states are instead given by phase variables~$\theta_k,\theta_j,\theta_l\in \mathbb{T}$, then nonpairwise interactions may have the form~$g(\theta_k,\theta_j,\theta_l)=\sin(\theta_k+\theta_j-2\theta_l)$. These interactions are characteristic of nonpairwise extensions of the Kuramoto model~\cite{ashwin2016hopf,leon2025theory}. Thus, what constitutes a nonpairwise interaction depends on the description of the system.

Interactions at the level of physical coordinates and nonpairwise interactions at the phase level are related through phase reduction for weakly coupled systems. Phase reduction is a technique that reduces a high-dimensional oscillatory system to a lower-dimensional phase description by exploiting the existence of stable limit cycles~\cite{nakao16,monga2019phase,pietras2019network,kuramoto2019concept}. Instead of tracking the full state of each oscillator, one describes its evolution in terms of a single phase variable. Phase reductions are computed as expansions in the coupling strength~$\varepsilon$, which quantifies how strongly oscillators interact. Classical approaches, including Kuramoto's original work~\cite{kuramoto1975}, are first order approximations in~$\varepsilon$ and yield phase models with pairwise interactions, such as the Kuramoto or Winfree models~\cite{winfree1967biological,Kuramoto_book}. These models have been widely used to study synchronization transitions and collective behavior~\cite{acebron2005kuramoto,rodrigues2016kuramoto}.

More recently, different approaches have been developed to compute higher-order phase reductions beyond first order~\cite{ashwin2016hopf,leon19,nijholt2022emergent,gengel2020high,bick2024higher,mau2024phase,fujii2026emergence}. These approaches systematically compute corrections to the phase dynamics and show that additional interaction terms appear at higher orders in~$\varepsilon$. In particular, nonpairwise phase interactions emerge naturally at second order and beyond, even when the original physical interactions are pairwise~\cite{leon19,bick2024higher,mau2024phase}. Therefore, a nonpairwise phase interaction does not necessarily imply the presence of a nonpairwise interaction in the original physical system. At the same time, if the physical interactions are themselves nonpairwise, these can already generate nonpairwise terms at first order in the phase reduction~\cite{leon2024,leon2025theory,leon2026symmetry}.

There are therefore two different routes through which nonpairwise phase interactions can arise. First, nonpairwise phase interactions can arise directly from nonpairwise interactions in the physical coordinates. If the physical interactions are nonpairwise, one generally expects nonpairwise interactions already in the first order phase dynamics. We refer to these as \textit{physical nonpairwise (PN) interactions} . Second, even if the physical interactions are pairwise, nonpairwise phase interactions can emerge at second and higher orders in the phase equations. These interactions arise from the reduction procedure itself and from the nonlinear response of the oscillators to the physical coupling. We refer to these as \textit{emergent nonpairwise (EN) interactions}.

Often these two routes have been considered separately. Some studies start from physical systems with explicit nonpairwise coupling and derive the corresponding phase equations~\cite{leon2024,leon2025theory,leon2026symmetry}. Other studies assume pairwise physical interactions and focus on nonpairwise terms emerging from higher-order phase reduction~\cite{leon19,bick2024higher}. In addition, nonpairwise phase interactions have often been introduced directly at the phase level as models, without linking them explicitly to an underlying physical system or reduction procedure~\cite{tanaka2011multistable,bick2016chaos,skardal2019abrupt,skardal2020higher,millan2020explosive,lucas2020multiorder,kovalenko2021contrarians,adhikari2023synchronization,skardal2023multistability,carballosa2023cluster,costa2024bifurcations,huh2024critical,wang2024coexistence,smith2024determining,dai2025higher}. As a result, the relationship between PN and EN remains only partially understood. In particular, it is natural to ask whether the same nonpairwise phase harmonic can originate from the two mechanisms, whether its physical origin can be identified from the resulting phase dynamics, and whether the two types of interactions can be used together.

To address these questions, here we exploit a recently developed parametrization method for higher-order phase reduction~\cite{von2023parametrisation}. Rather than deriving the phase dynamics alone, the method computes an approximation of the invariant torus of the coupled system together with the dynamics on it. The phase reduction and the embedding of the perturbed invariant torus are obtained order by order from the same conjugacy equation, i.e., the equation that maps the oscillator state onto its phase description. In this way, the quantities computed at a given order provide the inhomogeneity needed at the following order. This makes the method particularly suitable for computing higher-order corrections and, in the present context, for keeping track of the different contributions that generate nonpairwise phase interactions.

We first revisit the parametrization method and relate it to classical isochrone-based phase reduction. We then apply it to Stuart--Landau oscillators to compute the phase reduction up to second order. This allows us to consider PN and EN within the same framework and to compare their phase representations explicitly. We show that the two mechanisms can generate analogous nonpairwise harmonics, while differences in their complete expressions allow their origin to be distinguished in the system considered here.

Finally, we use these insights to develop a mixed order phase reduction that can be used to design effective phase coupling. 
A physical oscillator system with pairwise coupling generally does not behave exactly as its first order phase approximation at finite coupling strength: Second order terms, including EN, modify the phase dynamics, which is particularly relevant when the first order reduction is degenerate. 
We therefore ask whether PN can be introduced in the physical system to compensate for some of these corrections. 
For the Stuart--Landau system considered here, we use the correspondence between physical nonpairwise couplings and their phase interactions to engineer PN that oppose specific EN generated at second order. In this way, we modify the original physical system so that its dynamics are closer to those of its first order Kuramoto approximation.

\section*{Results}
\label{sec:results}

\subsection*{Higher-order phase reductions via parametrization}

Parameterization of invariant manifolds provides a geometric way to compute phase reductions systematically~\cite{von2023parametrisation,Bick2025}.
Here we compare this approach to traditional phase reduction techniques (see Methods for details and an example how to compute it for coupled Stuart--Landau oscillators).

Consider~$n$~identical uncoupled oscillators~$\dot{\vec{x}}_k = \vec{F}(\vec{x}_k)$, each one individually evolving on a stable limit cycle~$\Gamma\in\mathbb{R}^d$ in physical state space---we assume that they are identical for simplicity of notation.
Together, the dynamics of the~$n$~oscillators with joint state~$\vec{\mathbf{x}}=(\vec{x}_1, \dotsc, \vec{x}_n)$ evolve on an~$n$-dimensional torus~$\Gamma^n\subset\mathbb{R}^{nd}$ in the joint state space of all oscillators.
This torus can now be parameterized through phase variables---the natural choice of coordinates on a torus---through a map~$\vec{e}^{(0)}_k(\vartheta)$, which maps a phase~$\vartheta\in\mathbb{T} = \mathbb{R}/2\pi\mathbb{Z}$ to the corresponding point on~$\Gamma$ in physical space.
One has to make a choice on how to define the phase~$\vartheta$: One typically defines ~$\vec{e}^{(0)}_k(\vartheta) = \vec{\gamma}(2\pi\vartheta/T)$, where~$\vec{\gamma}(t)$~is the periodic solution with period~$T$, so that the uncoupled oscillator evolves at uniform speed~$\omega_k=2\pi/T$.
Together, these maps parameterize the invariant torus~$\vec{e}^{(0)}(\mathbb{T}^n) = \Gamma^n$.

The existence of an attracting invariant torus means that the dynamics can be reduced to phases: The dynamics near the torus are governed by the dynamics of the phases~$\vec{\theta}=(\theta_1, \dotsc, \theta_n)^\top$ only and ``amplitudes'' are fully determined by the phases.
For coupled oscillators 
\begin{align}\label{eq:Unred}
    \dot{\vec{x}}_k &= \vec{F}(\vec{x}_k) + \eps \vec{G}(\vec{\mathbf{x}}),
\intertext{
general mathematical theory~\cite{Fenichel1972} implies that there is an invariant torus~$\ord{\mathsf{T}}{\varepsilon}$ near~$\Gamma^n$ for nonzero coupling strength~$\varepsilon$.
Parameterization goes beyond existence: It allows us to compute not only the dynamics on the~$n$-dimensional torus (the phase reduction) given by} \dot{\vec{\theta}} &= \vec{f}(\vec{\theta})
\end{align} 
but also the relationship between phases and amplitudes through a map~$\vec{e}(\mathbb{T}^n) = \ord{\mathsf{T}}{\varepsilon}$. 
Reduced and unreduced dynamics are related through the conjugacy equation\footnote{Note that~$\vec{e}{~}'$ is a matrix that maps a vector into the tangent space at the torus.}
\begin{equation}\label{eq:ConjEqMain}
    \vec{e}{~}'\cdot \vec{f} = (\vec{F} + \eps \vec{G})\circ \vec{e},
\end{equation}
where the derivative~$'$ denotes~$\partial_{\vec{\theta}}$, which equates the dynamics on the invariant torus with the unreduced dynamics~\eqref{eq:Unred}; we omit the dependence of~$\vec{e}$ and~$\vec{f}$ on the phase vector~$\vec{\theta}$ to lighten the notation.
In practice, this equation can be solved by an expansion of the phase dynamics~$\vec{f} = \vec{\omega}+\varepsilon \vec{f}^{(1)}+\dotsb$ and the embedding~$\vec{e} = \vec{e}^{(0)} + \varepsilon \vec{e}^{(1)}+\dotsb$, substituting this into the conjugacy equation, collecting terms of equal order, and then solve the resulting system order-by-order.

When solving the conjugacy equations, one has to choose the phase variables along the torus.
If we split the embedding~$\vec{e}$ into a part~$\vec{g}$ tangential to the torus and a part~$\vec{h}$ normal to the torus (see Methods), this corresponds to a choice of~$\vec{g}$.
From the perspective of classical phase reduction, the standard choice is to let~$\vec{g} = \check{\vec{g}} =\ord{\vec{g}}{0}$ for any~$\eps$. In other words, it is fully determined by~$\ord{\vec{e}}{0}$, which is chosen so that~$\ord{\vec{f}}{0}=\vec{\omega} = (\omega_1, \dotsc, \omega_n)^\top$ is the uniform rotation as mentioned above.
Another possible choice is to parameterize in a way that removes nonresonant terms from the phase dynamics~$f^{(\varepsilon)}$ to finite order.
This leads to a reparameterized notion of phase through a choice of~$\vec{g} = \hat{\vec{g}} = \ord{\vec{g}}{0}+\eps\ord{\hat{\vec{g}}}{1}+\eps^2\ord{\hat{\vec{g}}}{2}+\dotsb$ such that the order-$\ell$ terms~$\ord{\hat{\vec{g}}}{\ell}$ are not necessarily zero.
While the first choice of~$\vec{g} =\check{\vec{g}}$ retains the original meaning of phase of the uncoupled oscillators, the second choice~$\vec{g}=\hat{\vec{g}}$ introduces a new definition of phase that simplifies the phase-reduced vector field~$\ord{\vec{f}}{\varepsilon}$.

Going from one parameterization to the other corresponds to a near-identity transformation, which removes nonresonant terms in the phase dynamics~$\ord{\vec{f}}{\varepsilon}$ to finite order---an averaging transformation.
More explicitly, note that as parameterization of the tangential direction, both~$\check{\vec{g}}$ and~$\hat{\vec{g}}$ are invertible as maps into the tangent space of the (approximation of) the invariant torus~$\mathsf{T}^{(\eps)}$.
This gives a near-identity transformation~$\vec{q} = \hat{\vec{g}}\circ\check{\vec{g}}^{-1}$ given by
\begin{align}\label{eq:NearIdentity}
     \vec{q}(\vec{\theta}) = (\hat{\vec{g}}\circ\check{\vec{g}}^{-1})(\vec{\theta}) 
     &= \vec{\theta} +\eps\ord{\hat{\vec{g}}}{1}((\ord{\vec{g}}{0})^{-1})+\eps^2\ord{\hat{\vec{g}}}{2}((\ord{\vec{g}}{0})^{-1}(\vec{\theta}))+\dotsb.
\end{align}
By our choice, the near-identity transformation removes nonresonant terms from the phase dynamics~$\ord{\vec{f}}{\varepsilon}$ up to finite order.
Thus, phase reduction through parameterization is an all-in-one approach to phase reduction: 
First, it provides a systematic way to compute phase-reduction to arbitrary order;
second, it provides an approximation of the phase dynamics~$\ord{\vec{f}}{\varepsilon}$ to the desired order as well as well as amplitude information for a given phase~$\vec{\theta}$ through~$\ord{\vec{e}}{\varepsilon}$;
third, it gives an explicit expression for the near-identity averaging transformation, which removes nonresonant terms from the phase dynamics.

By contrast, traditional approaches to phase reduction~\cite{nakao16,pietras2019network} first work in coordinates of the uncoupled oscillators and then perform averaging as a second step---without necessarily computing the near-identity phase transformation explicitly.
But since the object to compute is the same, there are naturally connections between the parameterization method (based on computing invariant manifolds) employed here and more classical approaches.
We discuss these connections in more detail in the Methods.

\subsection*{Comparison of physical and emergent nonpairwise interactions}

A key question regarding nonpairwise interactions in the Kuramoto model is whether there are intrinsic differences in the form of emergent and physical nonpairwise interactions. 
Recall that \textit{emergent} nonpairwise (EN) interactions are nonpairwise phase interactions that emerge in higher-order phase reduction from pairwise interactions in the unreduced equations.
Naturally, EN interactions arising at second order are an order of magnitude weaker than the pairwise ones~\cite{leon19,gengel2020high,bick2024higher}.
On the other hand, \textit{physical} nonpairwise (PN) interactions are a consequence of the structure of the interactions in the full model, which are already nonpairwise and, hence, present in the first order phase approximation~\cite{leon2024}. 

The actual relative strength of EN and PN in the phase equations, however, also depends on other factors, including intrinsic properties of the oscillators or the network coupling. 
For example, certain nonpairwise phase interactions to first order may vanish due to symmetries of the vector field and the coupling~\cite{leon2025theory,leon2026symmetry}.
Thus, we focus on the more general question whether EN and PN have distinct functional form.

To answer this question, we look at the smallest network configuration of three oscillators where nonpairwise interactions between all three oscillators appear.
More precisely, we focus on coupled Stuart--Landau (SL) oscillators with states given by~$z_j\in\mathbb{C}$ that evolve according to
\begin{equation}\label{eq:SL}
   \begin{split}
     \dot{z}_1=(a+ib)z_1+(c+id)\lvert z_1 \rvert^2z_1+\varepsilon G_1(z_1,z_2,z_3) , \\
     \dot{z}_2=(a+ib)z_2+(c+id)\lvert z_2 \rvert^2z_2+\varepsilon G_2(z_1,z_2,z_3), \\
     \dot{z}_3=(a+ib)z_3+(c+id)\lvert z_3 \rvert^2z_3+\varepsilon G_3(z_1,z_2,z_3),
 \end{split}
\end{equation}
for a general coupling function~$\vec G=(G_1, G_2, G_3)$, which determines the effect of~$z_1,z_2,z_3$ on a given oscillator.
The real parameters~$a, b, c, d,$ determine the intrinsic properties of each SL oscillator: 
For~$a>0$ and~$c<0$, Eq.~\eqref{eq:SL} admits a stable hyperbolic limit cycle~$\Gamma = \{\lvert z_k \rvert=R\}$~$z_0(t)=R\mathrm{e}^{i\omega t}$, with radius~$R = \sqrt{-\frac{a}{c}}$, frequency~$\omega=b-\frac{ad}{c}$, and Floquet exponent~$\lambda = -2a$.

In the following, we systematically compute phase interactions using the parameterization method outlined above.
This not only allows us to compare EN and PN phase interactions, but also yields a geometric interpretation of their origin compared to isochrone-based phase reduction approaches~\cite{leon19}.

\paragraph{Linear Coupling}
We first compute the emergent nonpairwise interactions that arise for linear coupling between the SL oscillators~\eqref{eq:SL} and an arbitrary network structure. 
To this end, consider pairwise physical network coupling of the form
\[G_j = G_j^\text{lin} := \mathrm{e}^{i\varrho}\sum_{k=1}^3w_{jk}z_k,\]
where~$w_{jk}\in\mathbb{R}$ is the (pairwise) interaction from oscillator~$k$ to~$j$ and~$\varrho$~is a rotational component to the coupling.
The first order phase interactions~$\vec f^{(1)}$ arising yield the Kuramoto-Sakaguchi model~\cite{sakaguchi1986soluble} with pairwise coupling terms that mirror the physical coupling (see Eq.~\eqref{eq:f1_expl_AppA} in Methods).
The second order phase interactions for the general system~\eqref{eq:Unred} are determined by the intrinsic dynamics~$\vec F$, the physical network coupling~$\vec G$, the first order phase dynamics~$\vec f^{(1)}$, the embedding~$\vec e$ of the torus up to first order, as well as their derivatives.
More precisely (see Eq.~\eqref{eq:homological2} in Methods), the second order phase interactions are determined by
\begin{align}\label{eq:InhomMain}
\vec U &= \vec G'(\ord{\vec e}{0}) \cdot  \ord{\vec e}{1}, &
\vec V &= -(\ord{\vec e}{1})' \cdot \ord{\vec f}{1}, &
\vec W &= \frac{1}{2}\vec F''(\ord{\vec e}{0})\left[\ord{\vec e}{1}, \ord{\vec e}{1}\right].
\end{align}
Since the computation in full generality becomes quite elaborate, we focus on the case of straight isochrones (i.e.,~$d=0$),~$c=-1$, and no self-coupling, i.e.,~$w_{jj}=0$.

The first source of EN interactions relates to the term~$\vec U$ and captures how the coupling function responds to the first order deformation of the torus:~$\vec G{}'(\ord{\vec e}{0})$ is the linearized coupling on the unperturbed torus (parameterized by~$\vec e^{(0)}$) that is applied to the first order deformation of the torus~$\ord{\vec e}{1}$ to measure the change in the coupling force induced by the torus already having been deformed at first order.
The resulting phase interaction terms for oscillator~$1$ are given by
\begin{equation}\label{eq:lin_second_C}
    \frac{1}{4a}\Big[w_{12}w_{21}\mathcal{C}_{1221}+w_{12}w_{23}\mathcal{C}_{1223} +w_{13}w_{31}\mathcal{C}_{1331}+w_{13}w_{32}\mathcal{C}_{1332} \Big],
\end{equation} 
where
\begin{displaymath}
    \begin{split}       \mathcal{C}_{1221}&=\sin(2\varrho)+\sin\Big(2(\theta_2-\theta_1)\Big) ,\\
       \mathcal{C}_{1223}&=\sin(\theta_3-\theta_1+2\varrho)+\sin(2\theta_2-\theta_3-\theta_1),\\
\mathcal{C}_{1331}&=\sin(2\varrho)+\sin\Big(2(\theta_3-\theta_1)\Big) , \\
       \mathcal{C}_{1332}&=\sin(\theta_2-\theta_1+2\varrho)+\sin(2\theta_3-\theta_2-\theta_1).
    \end{split}
\end{displaymath} 

\begin{figure}
    \centering
    \includegraphics[width=1\linewidth]{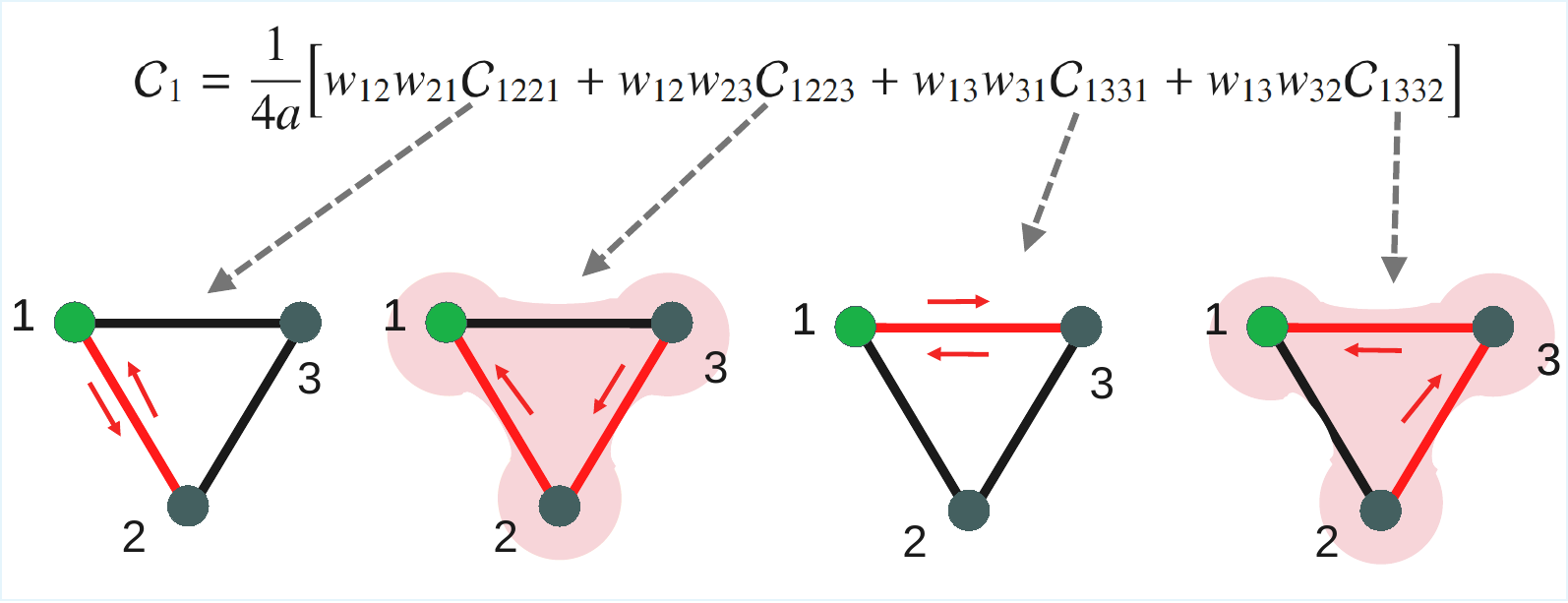}
    \caption{Pictorial representation of the emergent pairwise and nonpairwise interaction motifs from the contribution~$\vec{U}$ of the second order phase reduction.}
    \label{fig:inhomo1}
\end{figure}

The resulting phase interactions correspond to network motifs that are paths of length two; see Fig.~\ref{fig:inhomo1} for a pictorial representation.
Note that~$\mathcal{C}_{1223}$ (and, \textit{mutatis mutandis},~$\mathcal{C}_{1332}$) is an emergent nonpairwise (EN) interaction of oscillator~$3$ (resp.~$2$) towards oscillator~$1$ via oscillator~$2$ (resp.~$3$). While this is obvious for the harmonic~$\sin(2\theta_2-\theta_3-\theta_1)$ (also known as~$(2,-1,-1)$~\cite{namura2026optimal} or \textit{asymmetric}\footnote{Note that the term \textit{asymmetric} refers to the asymmetry with respect to the permutation of~$2$ and~$3$, not on the nature of the interactions. In this work, all interactions are symmetric, meaning that the effect of oscillator~$j$ on oscillator~$k$ is the same as the effect of oscillator~$k$ on oscillator~$j$. Using the terminology of network science, we are dealing with dynamics on undirected networks and hypergraphs.}~\cite{battiston2026collective} interaction), we can understand the nonpairwise nature of the harmonic~$\sin(\theta_3-\theta_1+2\varrho)$ thanks to the weights. 
On the other hand, the term ~$\mathcal{C}_{1221}$ (and, \textit{mutatis mutandis},~$\mathcal{C}_{1331}$) is an emergent pairwise interaction, but it is an indirect one: In fact, it is the interaction of oscillator~$1$ with itself via oscillator~$2$ (resp.~$3$).

The second source of EN interactions come from the term~$\vec{V}$, which relates the first-order phase dynamics and the first-order deformation of the tangent direction of the deformed torus.
The first order embedding~$\ord{\vec{e}}{1}(\vec\theta)$ captures the deformed torus, and its derivative sends the first order phase dynamics~$\ord{\vec{f}}{1}$ to the tangent space of the perturbed torus.
Concretely, the resulting phase interactions
for oscillator~$1$ are given by
\begin{equation}\label{eq:lin_second_D}
   - \frac{1}{4a}\Big[w_{12}^2\mathcal{D}_{1212}+2w_{12}w_{13}\mathcal{D}_{1213} +w_{13}^2\mathcal{C}_{1313} \Big],
\end{equation} where
\begin{displaymath}
    \begin{split}
       \mathcal{D}_{1212}&=\sin\Big(2(\theta_2-\theta_1)+2\varrho \Big),\\
       \mathcal{D}_{1213}&=\sin(\theta_2+\theta_3-2\theta_1+2\varrho),\\
       \mathcal{D}_{1313}&=\sin\Big(2(\theta_3-\theta_1)+2\varrho \Big).
    \end{split}
\end{displaymath}

\begin{figure}
    \centering
    \includegraphics[width=0.9\linewidth]{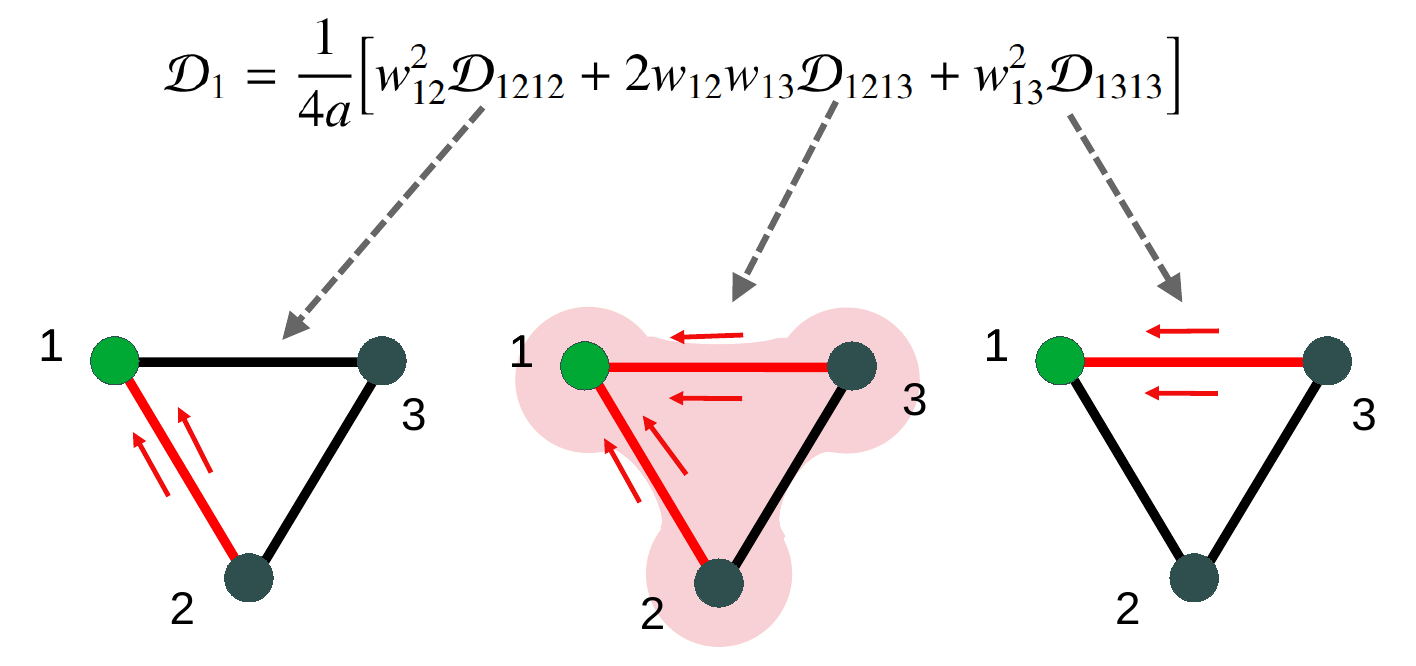}
    \caption{Pictorial representation of the emergent pairwise and nonpairwise interaction motifs from the contribution~$\vec{V}$ of the second order phase reduction.}
    \label{fig:inhomo2}
\end{figure}

The resulting phase interactions correspond to network motifs of length two that stay within a neighborhood of a given node; see Fig.~\ref{fig:inhomo2} for a pictorial representation.
The term~$\mathcal{D}_{1213}$ is an emergent nonpairwise (EN) interaction, namely, the simultaneous interaction of oscillators~$2$ and~$3$ towards oscillator~$1$ (also known as~$(1,1,-2)$~\cite{namura2026optimal} or \textit{symmetric}~\cite{battiston2026collective} interaction). On the other hand, the term ~$\mathcal{D}_{1212}$ (and, \textit{mutatis mutandis},~$\mathcal{D}_{1313}$), which stands for the interaction of oscillator~$2$ (resp.~$3$) on oscillator~$1$ taken two times, is an emergent pairwise interaction.
Note that the motifs correspond to these in Fig.~4 of Ref.~\cite{bick2024higher}, except for the one in panel~$(b1)$ as it actually cancels due to symmetry.

The potential source of EN which relates to the term~$\vec{W}$ is a curvature term intrinsic to the uncoupled dynamics.
Since~$\vec F$ is nonlinear, displacing a trajectory off the limit cycle by~$\vec e^{(1)}$ does not just translate the vector field but bends it, and~$\vec F''(\vec e^{(0)})$ quantifies this bending. 
For the case of straight isochrones, i.e.,~$d=0$, this contribution vanishes. 
This is because, for straight isochrones, the contribution generated by the second derivative of the uncoupled Stuart--Landau vector field is normal to the limit cycle and therefore does not contribute directly to the phase dynamics. 

The resulting second order phase reduction now combines pairwise first order terms with pairwise second-harmonic interactions and EN phase interactions at second order.
The first component of the second order phase vector field~$f^{(2)}
_1 = \big[
\vec f^{(2)}
\big]_1$ is
\begin{equation}
\label{eq:explicit_f2_first_component}
\begin{split}
f^{(2)}
_1
={}&
\frac{1}{4a}
\Big[
w_{12}w_{21}C_{1221}
+w_{12}w_{23}C_{1223}
+w_{13}w_{31}C_{1331}
+w_{13}w_{32}C_{1332}\\
&\hspace{15mm}
-w_{12}^{2}D_{1212}
-2w_{12}w_{13}D_{1213}
-w_{13}^{2}D_{1313}
\Big].
\end{split}
\end{equation}
This results in phase dynamics of oscillator~1 up to second order in the coupling strength~$\varepsilon$ given by
\begin{equation}
\dot{\theta}_1
=\omega
+\varepsilon
\Big(
w_{12}\sin(\theta_2-\theta_1+\varrho)
+w_{13}\sin(\theta_3-\theta_1+\varrho)
\Big)+\varepsilon^2 f^{(2)}
_1.
\label{eq:explicit_phase_dynamics_first_component}
\end{equation}
The equations for other oscillators are analogous.
Note that the EN phase interaction terms agree with what one would expect from a second order phase reduction~\cite{leon19} based on isochrones but without further geometric interpretation.

\paragraph{Nonlinear Coupling}

Next we consider how nonlinear physical couplings between Stuart--Landau oscillators~\eqref{eq:SL} are translated into PN phase interactions at first order.
The nonlinear interactions have to be at least cubic to be resonant (see Methods).
Again, we focus on the case of straight isochrones.

The first type of nonlinear interactions are given by the nonlinear coupling function
\begin{equation}\label{eq:cubic1}
G_j = G_j^\text{nln} := \mathrm{e}^{i\xi}\sum_{k,l=1\wedge k,l\neq j}^3 w_{jklk}z_k^2\bar{z}_l,
\end{equation}
where the network parameter~$w_{jkpq}$ corresponds to the (nonpairwise) interaction of the relevant nonlinear term and~$\xi$~is a rotational component to the coupling.

Applying the parametrization method yields the reduced first order phase dynamics
\begin{equation}\label{eq:reduced_cubic1}
   \begin{split}
     \dot{\theta}_1&=\omega+\varepsilon \Big(w_{1232}\sin(2\theta_2-\theta_3-\theta_1+\xi)+w_{1323}\sin(2\theta_3-\theta_2-\theta_1+\xi)\Big), \\
     \dot{\theta}_2&=\omega+\varepsilon \Big(w_{2131}\sin(2\theta_1-\theta_3-\theta_2+\xi)+w_{2313}\sin(2\theta_3-\theta_1-\theta_2+\xi)\Big), \\
     \dot{\theta}_3&=\omega+\varepsilon \Big(w_{3121}\sin(2\theta_1-\theta_2-\theta_3+\xi)+w_{3212}\sin(2\theta_2-\theta_1-\theta_3+\xi)\Big).
 \end{split}
\end{equation} 
Note that such nonpairwise phase interactions have been considered for example by Skardal and Arenas \cite{skardal2020higher} as they are amenable to exact mean-field reductions~\cite{Bick2018c}.

The second type of nonlinear interactions for the SL network~\eqref{eq:SL} are given by the nonlinear coupling function
\begin{equation}\label{eq:cubic2}
G_j = G_j^\text{nln} := \mathrm{e}^{i\xi}\sum_{k,l=1\wedge k,l\neq j}^3 w_{jkjl}z_k z_l\bar{z}_j,
\end{equation}

where we assume~$w_{jkjl}=w_{jljk}$ due to the symmetry of the interactions.
The first order phase reduction yields the following phase dynamics 
\begin{equation}\label{eq:reduced_cubic2}
   \begin{split}
     \dot{\theta}_1&=\omega+2\varepsilon w_{1213}\sin(\theta_2+\theta_3-2\theta_1+\xi), \\
     \dot{\theta}_2&=\omega+2\varepsilon  w_{2123}\sin(\theta_1+\theta_3-2\theta_2+\xi), \\
     \dot{\theta}_3&=\omega+2\varepsilon  w_{3132} \sin(\theta_1+\theta_2-2\theta_3+\xi),
 \end{split}
\end{equation} 
with PN phase interactions analogous to the nonpairwise harmonic considered by Tanaka and Aoyagi~\cite{tanaka2011multistable}.

\paragraph{Qualitative differences between EN and PN}

Both EN and PN yield the possible low-harmonic triplet interactions terms, the~$(1,1,-2)$ (``symmetric") and the~$(2,-1,-1)$ (``asymmetric") interactions.
But the results above allow us to identify key differences in the functional form of the EN and PN phase interactions for coupled SL oscillators.
First, the EN interactions are more constrained as relevant parameters such as phase lag\footnote{Note that we have denoted with~$\varrho$ and~$\xi$ the rotational components of the linear (pairwise) and nonlinear couplings, respectively, to distinguish them, but this would apply also if we choose the same rotational components for both couplings.} are functions of the linear (pairwise) physical coupling parameters.
Second, for EN, one can typically not have one type of nonpairwise interaction without the other type. As a concrete example, symmetric and asymmetric interactions arise to second order \emph{at the same time}.
Thus, generically one cannot expect exact mean-field reductions---such as Ott--Antonsen~\cite{ott2008low}---to work at higher-order.

In summary, one can generally expect that nonpairwise interactions in phase oscillator systems have contributions from both PN and EN phase interactions.
This overlap in the functional form of EN and PN provides an opportunity for designing phase coupling:
EN phase interactions could enhance or reduce PN interactions and vice versa.
We will explore this in the following section.

\subsection*{Mixed order phase reduction and coupling design}

In this Section, we exploit the results on second order phase reduction and the comparison between EN and PN interactions for two purposes. 
First, we test a mixed order phase reduction for systems with PN interactions. 
Indeed, a second order reduction for PN interactions is computationally challenging.
This raises the question of how well the original system can be approximated by retaining a second order reduction for linear pairwise interactions while stopping at first order for the nonpairwise ones. 
Second, we build on the insight gained from higher-order network motifs to pursue a higher-order coupling design approach, i.e., designing PN interactions so that the phase reduced equations take on a specific form.

\paragraph{Mixed-order phase reduction}

To understand the joint effect of EN and PN phase interactions, we focus on the network of~$N=3$ coupled Stuart--Landau oscillators~\eqref{eq:SL} without self-coupling ($w_{jj}=0$) and with straight isochrones ($d=0$).
We consider network interactions that decompose into a linear pairwise and a nonlinear nonpairwise part, each modulated by an independent small parameter.
Concretely, for small parameters~$\eps,\eta$ write
\begin{equation}
    G=\varepsilon G^{\mathrm{lin}}+\eta  G^{\mathrm{nln}}.
\end{equation} 
The idea is that the linear pairwise interactions~$G^{\mathrm{lin}}$ generate EN phase interactions to second order in~$\eps$, while~$G^{\mathrm{nln}}$ generate PN phase interactions to first order in~$\eta$.
Note that a second-order phase reduction for nonlinear coupling becomes computationally challenging even for SL oscillators with straight isochrones.

This provides the motivation for a mixed-order phase reduction approach: Consider the specific example that combines linear coupling with ``asymmetric'' (i.e.,~$(2,-1,-1)$) nonpairwise coupling discussed above, so that the explicit expression of the coupling for oscillator~$j=1$ is 
\begin{align*}   
G_1(z_1, z_2, z_3) = \varepsilon\mathrm{e}^{i\varrho}\sum_{k=2}^3w_{1k}z_k + \eta \mathrm{e}^{i\xi}\sum_{k,l=2}^3 w_{1klk}z_k^2\bar{z}_l.
\end{align*}
The mixed-order phase reduction (or simply ``$1.5$'') is now given by combining a second-order phase reduction (in~$\eps$) for linear coupling ($\eta=0$) with a first-order phase reduction (in~$\eta$) for nonlinear coupling ($\eps=0$).
Combining the expressions given in the previous section (see the Methods for a detailed derivation) yields the phase reduced dynamics for oscillator~$j=1$ as
\begin{align}\label{eq:1.5reduction}
    \dot{\theta}_1 &= 
    \omega+\varepsilon
\Big(
w_{12}\sin(\theta_2-\theta_1+\varrho)
+w_{13}\sin(\theta_3-\theta_1+\varrho)
\Big)+ \varepsilon^2 
f^{(2)}_1  \nonumber \\ &  
\quad\  +\eta \Big(w_{1232}\sin(2\theta_2-\theta_3-\theta_1+\xi)+w_{1323}\sin(2\theta_3-\theta_2-\theta_1+\xi)\Big), 
\end{align}
with~$f^{(2)}_1$ given by Eq.~\eqref{eq:explicit_f2_first_component}.

\begin{figure}
    \centering
   \includegraphics[width=1\textwidth]{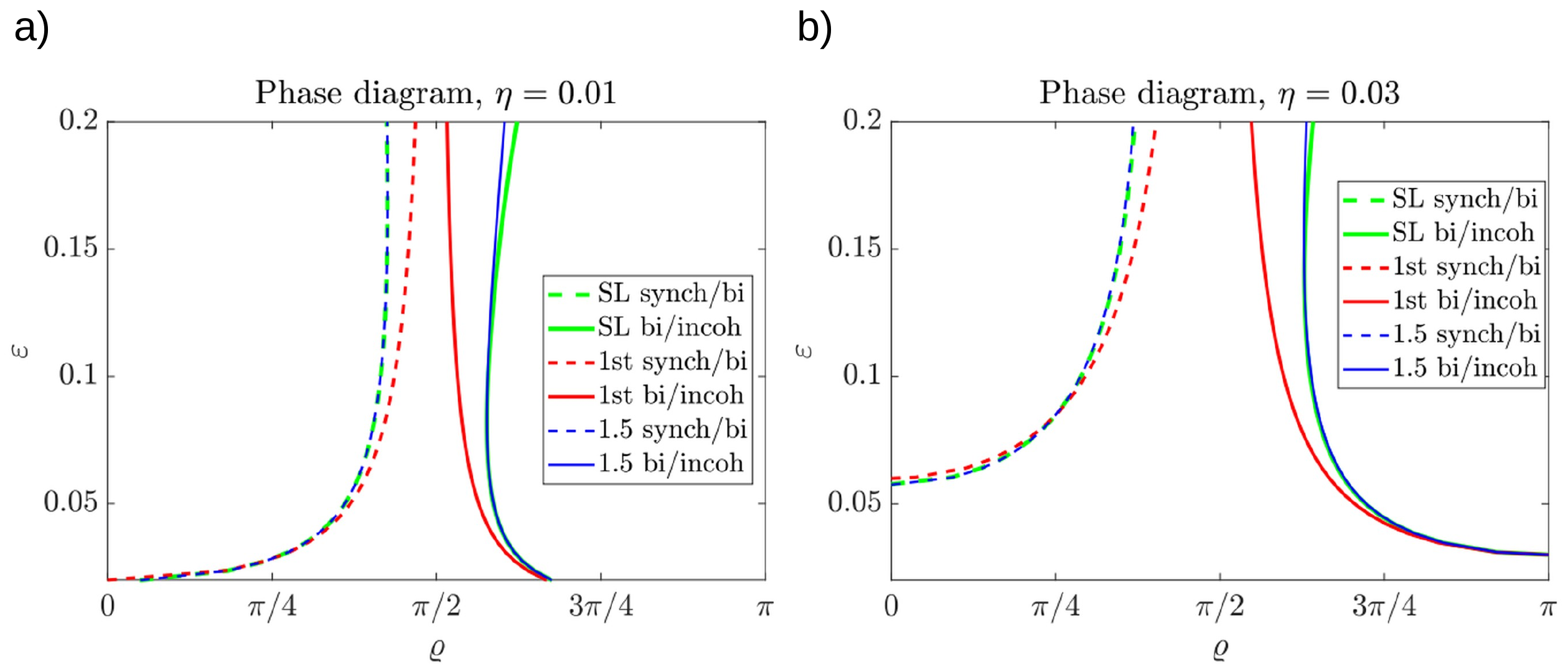}
    \caption{Comparison of the full Stuart--Landau system with linear and nonlinear coupling, its first order phase reduction, and the mixed ($1.5$) order phase reduction of Eq.~\eqref{eq:1.5reduction}, in the~$(\varrho,\varepsilon)$ plane, at fixed~$\xi=0$. Panels~$a)$ and~$b)$ show the boundaries separating the synchronized, bistable, and incoherent regions, obtained from the linear stability of the synchronized and twisted states (see Methods), for~$\eta=0.01$ and~$\eta=0.03$, respectively. Dashed lines mark where the twisted state becomes stable (boundary between the synchronized and bistable regions), solid lines where the synchronized state loses stability (boundary between the bistable and incoherent regions), and colors denote the model (green: full system; red: first order reduction; blue: mixed order reduction). The parameters are~$a=1$,~$b=1$,~$c=-1$,~$d=0$, all~$w$ are set to 1 except for the self couplings, which are set to 0.}
    \label{fig:mixed_rho_eps}
\end{figure}

\begin{figure}
    \centering
      \includegraphics[width=1\textwidth]{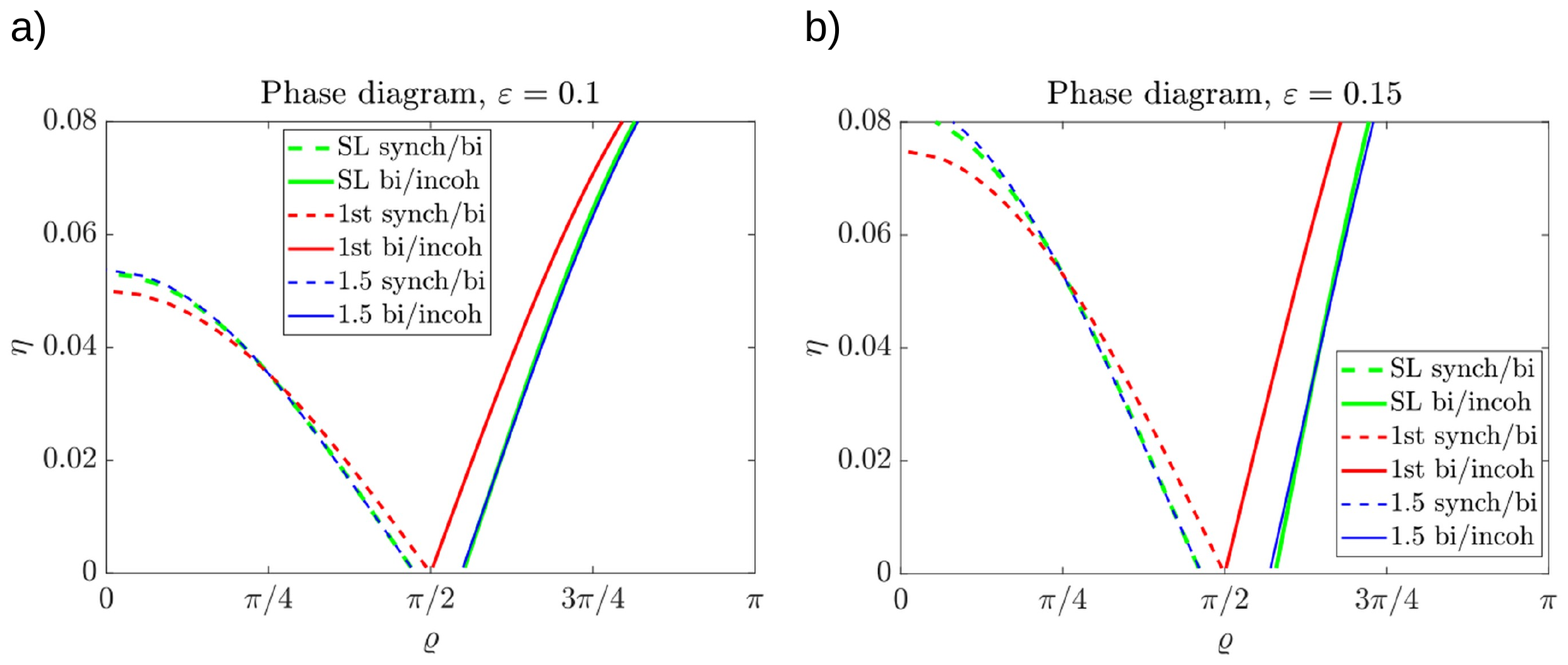}    
      \caption{Comparison of the full Stuart--Landau system with linear and nonlinear coupling, its first order phase reduction, and the mixed ($1.5$) order phase reduction of Eq.~\eqref{eq:1.5reduction}, in the~$(\varrho,\eta)$ plane, at fixed~$\xi=0$. Panels~$a)$ and~$b)$ show the boundaries separating the synchronized, bistable, and incoherent regions, obtained from the linear stability of the synchronized and twisted states (see Methods), for~$\varepsilon=0.1$ and~$\varepsilon=0.15$, respectively. Dashed lines mark where the twisted state becomes stable (boundary between the synchronized and bistable regions), solid lines where the synchronized state loses stability (boundary between the bistable and incoherent regions), and colors denote the model (green: full system; red: first order reduction; blue: mixed order reduction). The parameters are~$a=1$,~$b=1$,~$c=-1$,~$d=0$, all~$w$ are set to 1 except for the self couplings, which are set to 0.}
    \label{fig:mixed_rho_eta}
\end{figure}

The mixed-order phase reduction captures effects of both EN and PN phase interactions in the synchronization dynamics of the SL oscillators.
We compare the unreduced dynamics~\eqref{eq:SL} to the mixed-order phase reduction~\eqref{eq:1.5reduction} and a first order reduction (in both small parameters~$\eps,\eta$) without EN phase interactions.
Figure~\ref{fig:mixed_rho_eps} shows bifurcation boundaries as the strength~$\eps$ and rotational component~$\varrho$ is varied while the nonlinear coupling is fixed. 
We can observe that the mixed-order phase reduction gives a very good approximation of the full model. 
As expected, such approximation worsen as the coupling strength increases but it provides better approximation compared to only retaining first order terms.
Now, fixing the linear coupling strength~$\eps$,  Figure~\ref{fig:mixed_rho_eta} shows bifurcation boundaries as the nonlinear coupling~$\eta$ (and~$\varrho$ as above) is varied. 
Also here, we can see that the mixed-order gives good predictions of the bifurcation lines;
the approximation becomes worse as the coupling strengths increases.
Note that, for the mixed-order phase reduction to approximate the full system fairly well, we need the coupling of the PN to be about the order of the pairwise coupling squared, i.e.,~$O(\eta)\sim O(\varepsilon^2)$. 
Let us describe the phase diagrams in more detail. Since they are symmetric under~$\varrho\mapsto-\varrho$, it is enough to consider~$\varrho\in[0,\pi]$. The dashed and solid lines divide the diagrams into three regions: for small~$\varrho$ only the synchronized state is stable, for~$\varrho$ close to~$\pi$ only the twisted state is stable, and both states are stable in between (bistable region). Since the twisted state has $R=0$, we refer to the region where only the twisted state is stable as incoherent. By the symmetry, the synchronized region is centered at~$\varrho=0$ and the incoherent region at~$\varrho=\pi$. As the coupling strengths are varied, either of these regions may shrink to zero width. The corresponding boundary then reaches~$\varrho=0$ or~$\varrho=\pi$ and ends there. For example, in Fig.~\ref{fig:mixed_rho_eta}a) the dashed lines reach~$\varrho=0$ at~$\eta\approx0.05$. For larger~$\eta$, the twisted state is stable for every value of~$\varrho$. Hence, the synchronized region disappears, and the system is bistable for all~$\varrho$ to the left of the solid line.

\paragraph{Higher-order coupling design}

Finally, we can use the insights from higher-order phase reduction for \emph{coupling design}. 
Since nonpairwise phase interactions can arise both through linear coupling as EN interaction to higher order or through nonlinear coupling as PN interactions, we can use this interplay to enhance or suppress desired phase interaction terms in the phase reduced equations.
This is a step towards generalizing classical synchronization engineering~\cite{kiss2007engineering,kori2008synchronization,kiss2018synchronization} to higher orders and overcome a slow synchronization time scale induced by a first-order truncation.
Here we consider a more restricted design problem that combines higher-order phase reduction with design of nonpairwise phase interactions (see also~\cite{Liefting2023}) and will address the more general case in future work.

The concrete design problem is the following: As the strength~$\eps$ of linear coupling between SL oscillators is increased, which can allow for faster synchronization, the EN phase interactions become stronger.
Can we introduce weak nonlinear interaction so that the PN interactions cancel out the EN interactions?

We tackle this problem in the context of  three Stuart--Landau oscillators~\eqref{eq:SL} with straight isochrones,~$d=0$, and no self-coupling. 
With respect to oscillator~$1$, the EN phase interactions, which appear to second order in Eq.~\eqref{eq:explicit_f2_first_component}, are
proportional to~$\sin(2\theta_2-\theta_3-\theta_1)$,~$\sin(2\theta_3-\theta_2-\theta_1)$, and~$\sin(\theta_2+\theta_3-2\theta_1+2\varrho)$.
At the same time, nonlinear coupling of~$z_2^2\bar z_3$ yields PN phase interactions proportional to~$\sin(2\theta_2-\theta_3-\theta_1)$ to first order,~$z_3^2\bar z_2$ gives~$\sin(2\theta_3-\theta_2-\theta_1)$, and~$\mathrm{e}^{2i\xi}z_2z_3\bar z_1$ gives~$\sin(\theta_2+\theta_3-2\theta_1+2\xi)$. 
Since the functional form of the nonpairwise phase interactions overlap, the Fourier coefficients in the phase reduction will have contributions from both the linear coupling to second order as well as the nonlinear coupling to first order.

By choosing coefficients appropriately, nonlinear interactions can now be used to compensate for EN interactions arising to second order.
More specifically, consider
\begin{equation}
\begin{split}
\dot z_1 &=
(a+ib)z_1+c\lvert z_1\rvert^2 z_1
+\varepsilon \mathrm{e}^{i\varrho}
\left(w_{12}z_2+w_{13}z_3\right)
\\
&\qquad +\eta
\left[
-\left(
w_{1232}z_2^2\bar z_3
+w_{1323}z_3^2\bar z_2
\right)
+\mathrm{e}^{2i\varrho}w_{1213}z_2z_3\bar z_1
\right],
\end{split}\label{eq:engineered_SL_1}
\end{equation}
with analogous equations for oscillators~2 and~3. 
For EN and PN interactions to cancel in the corresponding mixed-order phase reduction, we choose
$w_{1232} = w_{12}w_{23}$,
$w_{1323} = w_{13}w_{32}$,
and
$w_{1213} = 2w_{12}w_{13}$,
to match the relative strength of EN and PN interactions
as well as matching phase-shift parameters, which we obtain by fixing~$\xi=\varrho$.

Given that the EN phase interactions appear at order~$\eps^2$ while the PN interactions are of order~$\eta$, the natural scaling relationship between the two small parameters is
\begin{equation}
    \eta=\frac{\varepsilon^2}{4a}.
    \label{eq:engineering_scaling}
\end{equation}
With this choice and the weights given above, the PN interactions cancel exactly both explicit EN harmonics in~$f^{(2)}_1$, i.e., the asymmetric~$(2,-1,-1)$ and the symmetric~$(1,1,-2)$ ones.
Corrections that mix the two physical couplings, such as~$O(\varepsilon\eta)$, enter only at~$O(\varepsilon^3)$ and are, therefore, beyond the second order approximation considered here: No second order phase reduction of the physical nonpairwise system is required to make the comparison consistent at this order.

The compensation of the second order corrections is nevertheless only partial, because PN and EN are not equivalent even when they generate the same harmonics.
The complete second order expression~\eqref{eq:explicit_f2_first_component} also contains the emergent pairwise interactions~$\mathcal{C}_{1221}$,~$\mathcal{C}_{1331}$,~$\mathcal{D}_{1212}$,~$\mathcal{D}_{1313}$, and the harmonics~$\sin(\theta_3-\theta_1+2\varrho)$ and~$\sin(\theta_2-\theta_1+2\varrho)$ in~$\mathcal{C}_{1223}$ and~$\mathcal{C}_{1332}$, whose nonpairwise nature is encoded in the weights only.
None of these terms can be compensated by nonpairwise couplings of the form~\eqref{eq:engineered_SL_1}.

\begin{figure}
    \centering
   \includegraphics[scale=0.28]{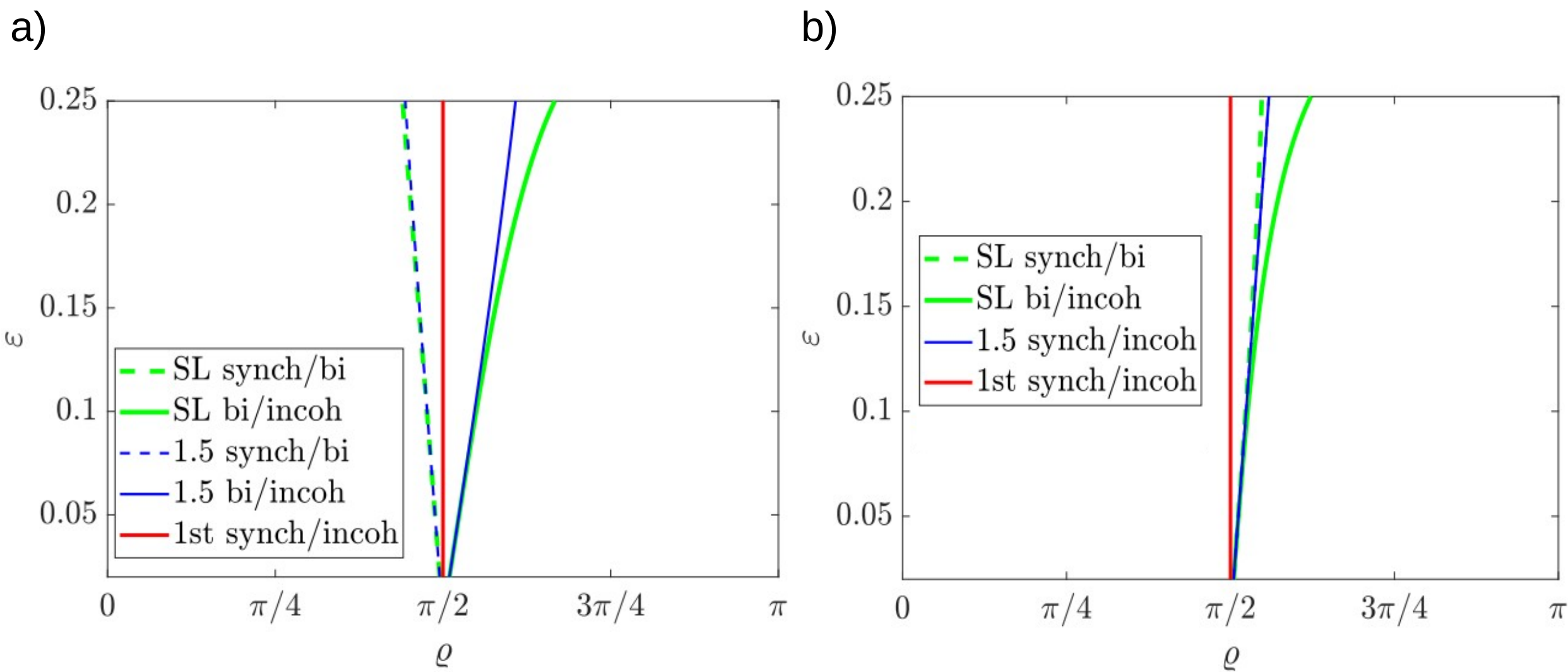}
    \caption{Phase diagram in the~$(\varrho,\varepsilon)$ plane for the Stuart--Landau model of Eq.~\eqref{eq:engineered_SL_1} (green), its first order phase reduction (red), and its mixed order phase reduction (blue), without $a)$ and with $b)$ the physical nonpairwise coupling. Dashed lines mark where the twisted state becomes stable (boundary between the synchronized and bistable regions); solid lines mark where the synchronized state loses stability (boundary between the bistable and incoherent regions), both obtained from linear stability (see Methods). The first-order phase reduction displays no bistability and therefore has a single boundary at~$\varrho=\pi/2$. $a)$~For~$\eta=0$, the mixed-order phase reduction, which coincides with the second-order for linear coupling, approximates the boundaries of the Stuart--Landau model. $b)$~For~$\eta=\varepsilon^2/(4a)$, Eq.~\eqref{eq:engineering_scaling}, with the weights chosen to cancel both explicit EN harmonics, the Stuart--Landau model has almost no bistability for small~$\varepsilon$, even though the curves do not coincide with the first order one. The parameters are~$a=1$,~$b=1$,~$c=-1$,~$d=0$; all~$w$ are set to~$1$ except for the self-couplings, which are set to~$0$. The mixed-order phase reduction displays no bistability, as the remaining second-order phase interactions are pairwise.}
    \label{fig:phase_diagram}
\end{figure}

To characterize the effect of the engineered PN interactions, we consider the phase diagram in the~$(\varrho,\varepsilon)$ plane, without ($\eta=0$, Fig.~\ref{fig:phase_diagram}a) and with ($\eta=\eps^2/(4a)$, Fig.~\ref{fig:phase_diagram}b) the engineered PN interactions. For the first order phase reduction, the transition between synchronization and incoherence occurs at~$\varrho=\pi/2$, independently of~$\varepsilon$, and no bistable region is present (red curves in Fig.~\ref{fig:phase_diagram}). For the original Stuart--Landau system, i.e., Eq.~\eqref{eq:engineered_SL_1} with~$\eta=0$, instead, finite coupling gives rise to a bistable region, delimited by the loss of stability of the synchronized and incoherent states (green curves in Fig.~\ref{fig:phase_diagram}a). As shown in Fig.~\ref{fig:phase_diagram}a, the mixed order phase reduction, which for~$\eta=0$ coincides with the second order one, reproduces the boundary between the synchronized and bistable regions of the non-reduced Stuart--Landau system with almost exactly, and the boundary between the bistable and incoherent regions for small and moderate~$\eps$ (blue curves in Fig.~\ref{fig:phase_diagram}a). The Stuart--Landau system of Eq.~\eqref{eq:engineered_SL_1} with~$\eta=\eps^2/(4a)$ (green curves in Fig.~\ref{fig:phase_diagram}b), instead, shifts both boundaries towards the first order Kuramoto transition: For small~$\eps$ the bistable region practically disappears. For larger~$\eps$ a narrower bistable region reappears and grows with~$\varepsilon$, due to corrections of order~$\eps^3$ that are not addressed by the design. The mixed order phase reduction (blue curves in Fig.~\ref{fig:phase_diagram}b), however, does not exhibit any bistability, as the nonpairwise interactions are designed to cancel out at second order and the remaining second order phase interactions are pairwise.
In other words, by designing the coupling appropriately we were able to suppress bistability even for larger coupling strengths~$\eps$.

\section*{Discussion}
\label{sec:discussion}

In this work, we have considered nonpairwise phase interactions arising from two different routes: Physical nonpairwise (PN) interactions already present in the physical coordinates, and emergent nonpairwise (EN) interactions generated by a second order phase reduction. First, we have revisited a recently introduced parametrization method to compute the phase reduction~\cite{von2023parametrisation}, and we have compared it with classical (isochrones-based) phase reduction. Then, we have put to use such method to compute a second order phase reduction and, hence, EN interactions, which we have compared with PN interactions. Lastly, we have exploited the knowledge of PN and EN interactions to design an oscillatory system that behaves as much as possible as a first-order phase model.

Our results show that PN and EN interactions cannot, in general, be identified only from the harmonics appearing in the phase dynamics. In fact, both mechanisms can generate the ``asymmetric''~$(2,-1,-1)$ and ``symmetric''~$(1,1,-2)$ nonpairwise harmonics. Nevertheless, their complete expressions are different. In the case of Stuart--Landau oscillators considered here, the phase shifts~$\varrho$ and~$\xi$ enter differently in PN and EN, and the second order reduction also generates terms whose nonpairwise nature is encoded in the coupling weights rather than explicitly in the phase variables. Hence, even when the same harmonic is present, its physical origin can be different. This distinction becomes less evident for particular choices of the parameters, for instance for real coupling, but the perturbative order and the dependence on the physical coupling still provide information about its origin.

From a methodological point of view, the parametrization method provides a direct way to compute these contributions order-by-order. In the present work, we have explicitly carried out the computation up to second order for linearly coupled Stuart--Landau oscillators. This already shows how nonpairwise terms arise from the different contributions to the second order inhomogeneity and how they can be related to the structure of the original pairwise interactions. At the same time, the explicit calculation shows that higher-order phase reductions rapidly become cumbersome, even for three Stuart--Landau oscillators. The possibility of organizing the calculation recursively through the conjugacy equation is therefore particularly useful when going beyond the first order.

For nonlinear couplings, however, a second-order phase reduction is considerably more computationally elaborate than for the pairwise case treated above: The inhomogeneity at second order would need to account for the nonpairwise structure of the coupling already present at first order. 
We have therefore not attempted this computation here. Instead, we propose a mixed order phase reduction obtained by retaining the second order terms for the linear pairwise coupling while stopping at first order for the nonpairwise one.
This balances the computational effort with a tractable approximation for systems with both types of interactions, without requiring the second-order reduction of the physical nonpairwise coupling.

The comparison between PN and EN interactions also suggests a possible use of physical nonpairwise interactions, namely, as a tool for coupling design. Since PN interactions enter already at first order in their own coupling strength, choosing their strength as~$\eta=O(\varepsilon^2)$ makes their contribution comparable with the EN generated at second order by the pairwise coupling.  For the Stuart--Landau example, we have chosen the physical nonpairwise couplings so that their phase reductions oppose the explicit~$(2,-1,-1)$ and~$(1,1,-2)$ EN interactions. The cancellation is necessarily partial, since PN and EN interactions are not equivalent and the second-order phase reduction contains additional terms that are not removed by the designed coupling. Nevertheless, this construction provides a simple example in which physical nonpairwise interactions can be used not to introduce additional collective dynamics, but rather to reduce the difference between the physical oscillator system and a prescribed first order phase model. The numerical results indicate that this compensation affects the synchronization region. In particular, while the second order phase reduction follows more closely the non-engineered Stuart--Landau system, the engineered physical system shifts towards the behavior predicted by the first order phase model. In this sense, the same distinction between PN and EN that prevents their exact identification also determines which second order effects can, and which cannot, be compensated through synchronization engineering.

In this work, we have performed the explicit computations of the second order phase reduction to identical Stuart--Landau oscillators with straight isochrones, without self-coupling and pairwise linear coupling. This setting allows us to distinguish the different contributions and to compare PN and EN interactions explicitly. However, let us remark that the parametrization method itself is not restricted to this simple setting: 
It provides a systematic way to compute higher-order phase reduction with a clear geometric interpretation and is related to Kuramoto's alternative phase reduction technique (see Methods).
Using parameterization to compute higher-order phase reduction does not require the specific functional form of the isochrones of an oscillator.
Hence, the main obstacle to treating more general systems is the amount of computation required rather than the method; this makes this framework well suited for future automation~\cite{Wedgwood2x}.
Moreover, we plan to exploit this method to compute higher-order phase reductions of systems with physical nonpairwise interactions, i.e., oscillator systems on hypergraphs, and understand the emerging motifs, as well as the case with delay, both with pairwise and (physical) nonpairwise interactions, so to extend the framework developed in~\cite{de2024higher,Bick2025,fujii2026emergence}. 
In the present work, PN interactions have been considered at first order, while EN interactions have been obtained from the second order reduction of pairwise coupling. 
Computing the second and higher orders when the physical system itself contains nonpairwise interactions would allow us to study how PN and EN interactions combine beyond the approximation considered here. 

Finally, the framework hereby developed can find further applications in synchronization engineering~\cite{kiss2007engineering,kori2008synchronization,rusin2009framework,rusin2010synchronization,kiss2018synchronization}, where different physical nonpairwise couplings could be designed to enhance or suppress specific terms emerging in the reduced phase dynamics. Nonpairwise couplings, although not explicitly referred to as such, have long been exploited in the design of oscillatory networks, for instance, to implement oscillatory neurocomputers and associative memories~\cite{hoppensteadt1999oscillatory,heger2016robust}. More recently, nonpairwise couplings have been explicitly implemented in electrical~\cite{vera2024electronic,minati2024chaotic} and electrochemical~\cite{nijholt2026hypernetworks} systems, which provides a promising avenue to bring theoretical insights to real-world systems.

\section*{Methods}
\label{sec:methods}

\subsection*{Higher-order phase reduction through parametrization}
\label{sec:Param}

Here we review the parameterization method to compute phase reductions introduced in~\cite{von2023parametrisation}.
Consider a network of coupled oscillators whose states~$\vec{x}_k\in\mathbb{R}^d$ evolve according to
\begin{align}\label{eq:Network}
    \dot{\vec{x}}_k &= \vec{F}(\vec{x}_k) + \eps \vec{G}(\vec{\mathbf{x}}),
\end{align}
where each isolated node is an oscillator, i.e.,~$\dot{\vec{x}}_k = \vec{F}(\vec{x}_k)$ has a hyperbolic limit cycle~$\Gamma$, and~$\vec{G}$~determines the interaction between nodes.
This means that the system~\eqref{eq:Network} for~$\eps=0$ has a normally hyperbolic torus~$\mathsf{T}^{(0)}=\Gamma^n$.
For sufficiently small coupling~$\eps>0$~\cite{Fenichel1972} this torus persists, that is, there is an invariant torus~$\mathsf{T}^{(\eps)}$ near~$\mathsf{T}^{(0)}$.

\paragraph{Example}
As a running example, we consider~$n$~coupled identical Stuart--Landau oscillators.
Identify~$\mathbb{C} \equiv \mathbb{R}^2$ and write~$z_k = x_k+iy_k \equiv (x_k,y_k)$ for the complex coordinates with~$i=\sqrt{-1}$ denoting the imaginary unit.
An (isolated) Stuart--Landau oscillator, whose state~$z_k\in\mathbb{C}$ evolves according to~$\dot z_k = F_k(z_k)$ with
\[F_k(z_k)=(a+ib)z_k+(c+id)\lvert z_k \rvert^2z_k,\]
and parameters~$a>0$ and~$c<0$ has a stable hyperbolic limit cycle~$\Gamma(t) = R\mathrm{e}^{i\omega t}$,
with radius~$R = \sqrt{-\frac{a}{c}}$, frequency~$\omega=b-\frac{ad}{c}$, and Floquet exponent~$\lambda = -2a$. The oscillators being identical, the frequency~$\omega$ is the same for each isolated system.
A network of~$n$ coupled Stuart--Landau oscillators now has a normally hyperbolic invariant torus~$\Gamma^n$ that persists for~$\varepsilon>0$.\\

The main idea of phase reduction through parameterization is to compute the perturbed invariant torus~$\mathsf{T}^{(\eps)}\subset\mathbb{R}^{m}$, where $m=nd$, and the phase reduction as the dynamics restricted to~$\mathsf{T}^{(\eps)}$.
Since natural coordinates for the~$n$-dimensional torus~$\mathsf{T}^{(\eps)}$ are given by~$\vec{\theta}=(\theta_1,\dotsc,\theta_n)^\top\in\mathbb{T}^n$, the phase reduction is of the form 
\begin{equation}\label{eq:PhaseRed}
    \dot\theta_k = f_k(\vec \theta).
\end{equation}
Thus, we need to compute a map~$\vec{e}:\mathbb{T}^n\to\mathbb{R}^m$ so that~$\vec{e}(\mathbb{T}^n) = \mathsf{T}^{(\eps)}$ and the phase dynamics~\eqref{eq:PhaseRed} determined by~$\vec{f}=(f_1, \dotsc, f_n)^\top$ match the dynamics of the unreduced system~\eqref{eq:Network} restricted to~$\mathsf{T}^{(\eps)}$.
In formulas, this is precisely the case if~$\vec{f},\vec{e}$ satisfy the \emph{conjugacy equation} (Eq.~\eqref{eq:ConjEqMain} in the Results)
\begin{equation}\label{eq:ConjEq}
    \vec{e}{~}'(\vec{\theta})\cdot \vec{f}(\vec{\theta}) = ((\vec{F} + \eps \vec{G})\circ \vec{e})(\vec{\theta}),
\end{equation}
which relates the phase dynamics~$\vec{f}$ (left hand side of the conjugacy equation) with the dynamics in the physical coordinates determined by~$\vec{F}, \vec{G}$ (right hand side) in~\eqref{eq:Network}. 
Note that the operator~$\vec e{~}'(\vec{\theta})$ sends the vector field on~$\mathbb{T}^n$ (the reduced system) to the tangent directions of~$\mathsf{T}^{(\eps)}$ in the physical coordinates, the derivative~$'$ denoting~$\partial_{\vec{\theta}}$.
Now expand~$\vec{f},\vec{e}$ in the small parameter~$\eps$,
\begin{equation} \label{eq:expansion_E_f}
\begin{aligned}
    \vec{f} &= \ord{\vec{f}}{0} + \eps \ord{\vec{f}}{1} +\dotsb + \eps^{\ell} \ord{\vec{f}}{\ell} + \dotsb \\
    \vec e & = \ord{\vec{e}}{0} + \eps \ord{\vec{e}}{1} +   \dotsb + \eps^\ell \ord{\vec{e}}{\ell} + \dotsb,
\end{aligned}
\end{equation} 
substitute into the conjugacy equations~\eqref{eq:ConjEq}, and Taylor-expand~$\vec{F},\vec{G}$.
For given order~$\ell\geq 0$, equating terms with~$\eps^\ell$ gives an inhomogeneous linear equation, the \emph{$\ell$th order homological equation}
\begin{align}\label{eq:lHom}
        \ord{\vec{e}}{\ell}(\vec{\theta})\cdot \ord{\vec{f}}{0}(\vec{\theta}) - \big(\vec{F}\big)'(\ord{\vec{e}}{0}(\vec{\theta})) [\ord{\vec{e}}{\ell}(\vec{\theta})] + \ord{\vec{e}}{0}{'}(\vec{\theta}) \cdot \ord{\vec{f}}{\ell}(\vec{\theta}) &= \ord{\vec{\eta}}{\ell}(\vec{\theta}).
\end{align}
The inhomogeneity~$\ord{\vec{\eta}}{\ell}(\vec{\theta})$ is a function of the network interactions, the embedding and phase dynamics, and their derivatives (e.g., the Jacobian matrices~$\big(\vec{F}\big)',\big(\vec{G}\big)'$ and Hessian tensor~$\vec{F}{}''$) from the multivariate Taylor expansion. 
Concretely, for~$\ell=1,2$ these are
\begin{align}
\ord{\vec{\eta}}{1}(\vec{\theta}) &= 
\vec{G}(\ord{\vec{e}}{0}(\vec{\theta})) \label{eq:homological1}\\
\ord{\vec{\eta}}{2}(\vec{\theta}) &= 
\big(\vec{G}\big)'(\ord{\vec{e}}{0}(\vec{\theta}))[\ord{\vec{e}}{1}(\vec{\theta})] + \frac{1}{2}\vec{F}{}''(\ord{\vec{e}}{0})[\ord{\vec{e}}{1}(\vec{\theta}), \ord{\vec{e}}{1}(\vec{\theta})] 
- \ord{\vec{e}}{1}{'}(\vec{\theta})\cdot\ord{\vec{f}}{1}(\vec{\theta}) \label{eq:homological2}
\end{align}
Importantly,~$\ord{\vec{\eta}}{\ell}(\vec{\theta})$ only depends on lower-order quantities, which means that they can be solved order-by-order as detailed in~\cite{von2023parametrisation}.

To solve the homological equations, it is natural to split the physical coordinates into directions tangential to the unperturbed invariant torus and normal to it.
The choice of coordinates in the tangential direction~$\mathbf{T}(\vec{\theta})=\ord{\vec{e}}{0}{'}(\vec{\theta}): \mathbb{R}^{n}\to\mathbb{R}^m$ directly relates how the torus is parameterized through~$\vec{\theta}$; indeed, choose~$\ord{\vec{e}}{0}(\vec{\theta})$ such that in the absence of coupling ($\eps=0$) the oscillators evolve at constant speed~$\dot{\theta}_k = \omega_k$, i.e.,~$\ord{\vec{f}}{0} = \vec{\omega} := (\omega_1, \dotsc, \omega_n)^\top$.
The choice of coordinates in the normal direction~$\mathbf{N}(\vec{\theta}):\mathbb{R}^{m-n}\to\mathbb{R}^m$ directly relates to the amplitude of the oscillation as it measures deviations from the unperturbed oscillators; the Floquet directions of the uncoupled limit cycles provide a natural choice of normal directions.

\begin{figure}
\centering\includegraphics[width=0.9\textwidth]{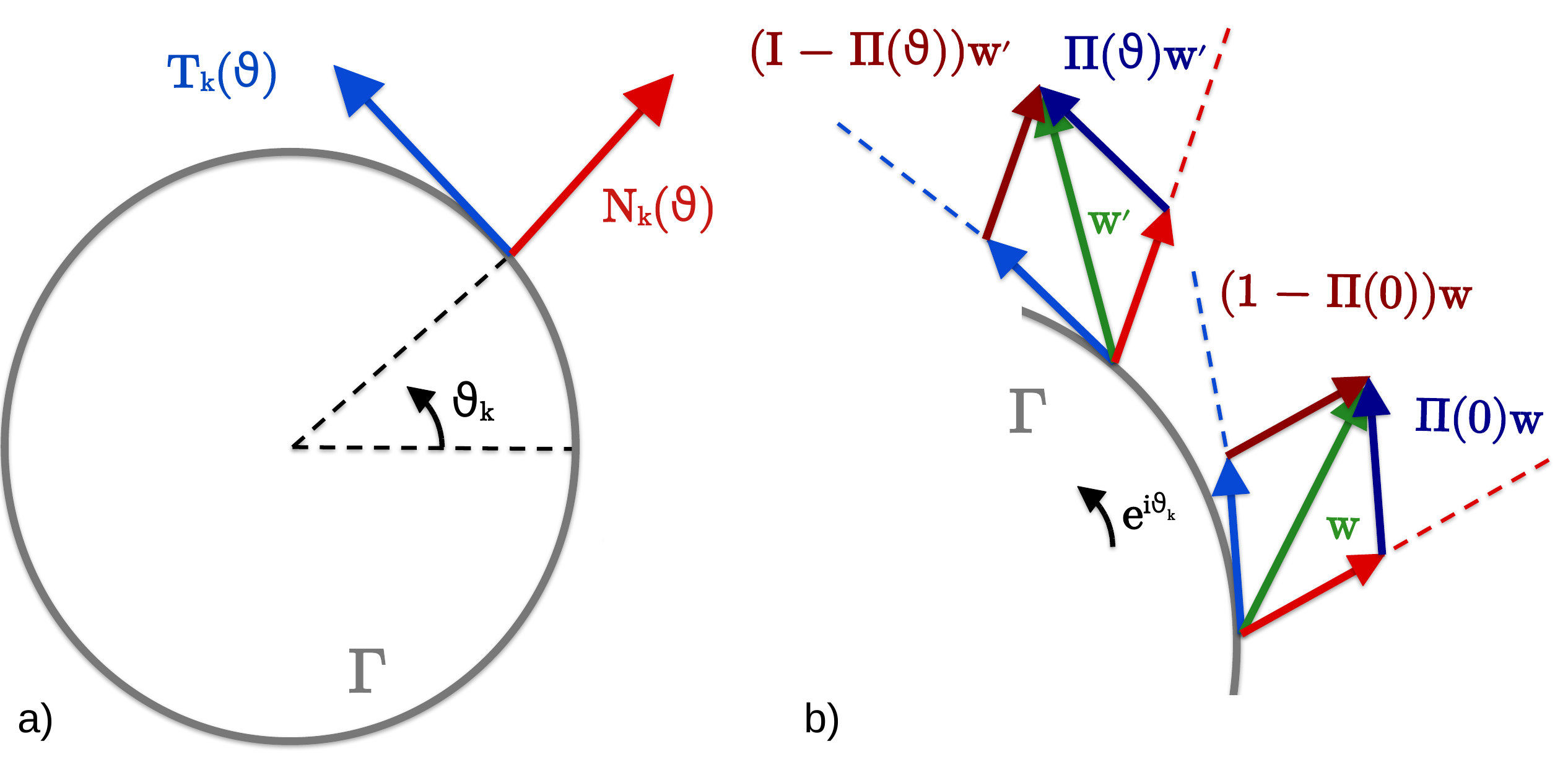}
\caption{Natural coordinates for a Stuart--Landau oscillator are given by tangential and normal directions. 
Panel~$a)$ shows the limit cycle~$\Gamma$ and the phase-dependent operators~$\mathbf{T}_k, \mathbf{N}_k$ that point in the tangential and normal directions, respectively.
Panel~$b)$ shows the projections~$\Pi, \mathbf{1}-\Pi$ on these directions. The projections~\eqref{eq:Pi0} for~$\theta_k=0$ are straightforward to compute.
Rotation by~$\theta_k$ corresponds to multiplication by~$\mathrm{e}^{i\theta_k}$; this yields~$\Pi, \mathbf{1}-\Pi$ by rotating by~$-\theta_k$, applying~$\Pi(0)$, and rotating back.
\label{fig:SLsketch}}
\end{figure}

\paragraph{Example}
Parameterizing the unperturbed limit cycle~$\Gamma$ of the Stuart--Landau oscillators in the first Example by
\[\ord{e}{0}_k = R\mathrm{e}^{i\theta_k}\]
yields the trivial zeroth order phase reduction with a uniformly rotating phase~$\dot\theta_k = \omega = 1$.

For an individual oscillator, the vector~$\mathbf{T}_k$ gives the tangent direction and~$\mathbf{N}_k$~gives a normal direction to the to the limit cycle~$\Gamma$ at a give phase; cf.~Fig.~\ref{fig:SLsketch}.
Concretely, at~$z_k = 1 \equiv (1, 0)$, the point of zero phase~$\theta_k = 0$, the tangent direction direction~$\mathbf{T}_k = i \equiv (0, 1)$ is vertical.
The normal direction~$\mathbf{N}_k = c+id \equiv (c,d)$ points into the direction of the isochrone; stability in this transverse direction is determined by the Floquet exponent~$\lambda = -2a = L_{kk}$.
Because of the rotational symmetry, the corresponding directions 
for a general phase~$\theta_k$ are obtained by rotation, that is,
\begin{align*}
\mathbf{T}_k &= \big(\ord{e}{0}_k\big)' = iR\mathrm{e}^{i\theta_k}, &
\mathbf{N}_k &= (c+id)\mathrm{e}^{i\theta_k}.
\intertext{
For the joint system of~$n$ oscillators, the corresponding matrices are diagonal}
\mathbf{T}(\vec\theta) &= \diag(iR\mathrm{e}^{i\theta_k}), &
\mathbf{N}(\vec\theta) &= \diag((c+id)\mathrm{e}^{i\theta_k}).
\end{align*}
and the Floquet matrix~$L(\theta) = \diag(-2a)$
is constant because of the rotational symmetry.\\

With the coordinates fixed, the embedding can now be split into a component~$\ord{\vec{g}}{\ell}: \mathbb{T}^n \to \mathbb{R}^{n}$ along~$\mathbf{T}(\vec{\theta})$ and a component~$\ord{\vec{h}}{\ell}: \mathbb{T}^n \to \mathbb{R}^{m-n}$ in the direction~$\mathbf{N}(\vec{\theta})$.
In other words, the embedding~$\ord{\vec{e}}{\ell}: \mathbb{T}^{n}\to\mathbb{R}^m$ can be written as
\[
\ord{\vec{e}}{\ell}(\vec{\theta}) = \mathbf{T}(\vec{\theta})\cdot\ord{\vec{g}}{\ell}(\vec{\theta}) + \mathbf{N}(\vec{\theta})\cdot\ord{\vec{h}}{\ell}(\vec{\theta}).
\]
With the Floquet matrix~$L$ capturing the dynamics in the normal directions, the~$\ell$th-order homological equation~\eqref{eq:lHom} in the new coordinates is
\begin{subequations}\label{eq:conj_order}
\begin{align}
    \mathbf{T}(\vec{\theta})\cdot
    \left[\partial_{\vec{\omega}} \ord{\vec{g}}{\ell}(\vec{\theta})  + \ord{\vec{f}}{\ell}(\vec{\theta})\right] &= \Pi(\vec{\theta})\cdot\ord{\vec{\eta}}{\ell}(\vec{\theta})\label{eq:ConjODETang} \\
    \mathbf{N}(\vec{\theta})\cdot\left[\partial_{\vec{\omega}} \ord{\vec{h}}{\ell}(\vec{\theta}) - L \ord{\vec{h}}{\ell}(\vec{\theta})\right] &= (\mathbf{1}-\Pi(\vec{\theta}))\cdot\ord{\vec{\eta}}{\ell}(\vec{\theta}), \label{eq:ConjODENorm}
\end{align}
\end{subequations}
where~$\partial_{\vec{\omega}}$ is the derivative in the direction~$\vec{\omega}$ and~$\Pi:\mathbb{R}^n\to\mathbb{R}^n$ and~$(\mathbf{1}-\Pi):\mathbb{R}^n\to\mathbb{R}^n$ are the projections in the tangential and normal directions, respectively.

\paragraph{Example}
For the Stuart--Landau oscillators in the first Example, the projections~$\Pi, \mathbf{1}-\Pi$ at phase~$\theta_k = 0$ ($z_k = 0$) of a vector~$w = u+iv \equiv (u,v)$ onto~$\mathbf{T}_k$ and~$\mathbf{N}_k$ are
\begin{subequations}
\label{eq:Pi0}
\begin{align}
\Pi_k(0)(u+iv) &=i\left(v-\left(\frac{d}{c}\right)u\right),\\
(1-\Pi_k(0))(u+iv) &=u+iu\left(\frac{d}{c}\right).
\end{align}
\end{subequations}
The rotational symmetry allows us to compute the corresponding projections by rotating~$w$ to zero phase, applying~$\Pi_k(0)$ and rotating back, that is,
\begin{align*}
\Pi_k(\theta_k)(u+iv) &= \mathrm{e}^{i\theta_k}\Pi_k(0)(u+iv)\mathrm{e}^{-i\theta_k},\\
(1-\Pi_k(\theta_k))(u+iv) &= \mathrm{e}^{i\theta_k}(1-\Pi_k(0))(u+iv)\mathrm{e}^{-i\theta_k}.
\end{align*}
For the system of~$n$~oscillators, the projections, which determine the right hand side of~\eqref{eq:conj_order}, are given by the matrices
\begin{align}\label{eq:SLProjections}
\Pi(\vec{\theta}) &= \diag(\mathrm{e}^{i\theta_k}\Pi_k(0)\mathrm{e}^{-i\theta_k})\\
\mathbf{1}-\Pi(\vec{\theta}) &= \diag(\mathrm{e}^{i\theta_k}(1-\Pi_k(0))\mathrm{e}^{-i\theta_k}).
\end{align}

Finally, the operators~$\mathbf{T}, \mathbf{N}$ can be inverted on their image to determine the embedding~$\ord{\vec{g}}{\ell}, \ord{\vec{h}}{\ell}$ and the phase dynamics~$\ord{\vec{f}}{\ell}$.
More precisely, with the pseudo-inverses~$\mathbf{T}^+, \mathbf{N}^+$, we obtain
\begin{subequations}\label{eq:ConjSolve}
\begin{align}
    \partial_{\vec{\omega}} \ord{\vec{g}}{\ell}(\vec{\theta})  + \ord{\vec{f}}{\ell}(\vec{\theta}) &= \mathbf{T}^+(\vec{\theta})\cdot \Pi(\vec{\theta})\cdot \ord{\vec{\eta}}{\ell}(\vec{\theta}), \label{eq:fSol}\\
    (\partial_{\vec{\omega}} - L )\ord{\vec{h}}{\ell}(\vec{\theta}) &= \mathbf{N}^+(\vec{\theta})\cdot (\mathbf{1}-\Pi(\vec{\theta}))\cdot \ord{\vec{\eta}}{\ell}(\vec{\theta}).
\end{align}
\end{subequations}
Thus, the components of the embedding~$\ord{\vec{g}}{\ell}, \ord{\vec{h}}{\ell}$ and the phase dynamics~$\ord{\vec{f}}{\ell}$ are determined by the quantities on the right-hand-side.
In other words, the choice of coordinates along the invariant torus determined by~$\ord{g}{\ell}$ determines the vector field~$\ord{\vec{f}}{\ell}$ and vice versa.
One option is to simply set~$\ord{\vec{g}}{\ell}=0$---this corresponds to phase coordinates in terms of the uncoupled limit cycles---which gives~$\ord{\vec{f}}{\ell}(\vec{\theta}) = \mathbf{T}^+(\vec{\theta})\cdot \Pi(\vec{\theta})\cdot \ord{\vec{\eta}}{\ell}(\vec{\theta})$.
Alternatively, one can choose~$\ord{\vec{g}}{\ell}$ so that~$\ord{\vec{f}}{\ell}$~contains only resonant combinations of phase angles terms (e.g.,~$\sin(\theta_j-\theta_k)$) and no nonresonant terms (e.g.,~$\sin(\theta_k)$); see~\cite{von2023parametrisation}.
This corresponds to a coordinate change of~$\vec{\theta}$ which simplifies the phase dynamics~$\ord{\vec{f}}{\ell}$: 
With this choice of coordinates, the phase reduced vector field~$\vec{f}$~is in \emph{normal form}.

\paragraph{Example}
For the Stuart--Landau oscillators in first Example, the pseudo inverses of~$\mathbf{T}, \mathbf{N}$ are given by
\begin{align}\label{eq:SLInverses}
\mathbf{T}^+(\vartheta) &= \diag((iR)^{-1}\mathrm{e}^{-i\theta_k}),&
\mathbf{N}^+(\vartheta) &= \diag((c+id)^{-1}\mathrm{e}^{-i\theta_k}).
\end{align}
This is all the information needed to compute first and second order phase approximations for specific linear and nonlinear coupling~$\vec G$ between oscillators.\\

\subsection*{Comparison of parameterization with existing phase reduction techniques}
\label{sec:CompPhaseRed}

\newcommand{\la}{\langle}
\newcommand{\ra}{\rangle}

\newcommand{\pa}{\partial}

In this section, we briefly discuss the relationship between the parametrization method and the classical phase reduction methods.
First, the simplest first order phase reduction method is due to Winfree~\cite{winfree1967biological} and Kuramoto~\cite{Kuramoto_book}. It is based on the notion of asymptotic phase and called Method~I in~\cite{Kuramoto_book}.
It is interesting to note that a more systematic higher-order phase reduction method is also developed in~\cite{Kuramoto_book}, which is called Method~II.
Here, we discuss the relation between the Kuramoto's method and the Parametrization Method discussed in this paper. The essential ideas for setting up the perturbation expansion are exactly the same. However, how to solve the resulting equations at each order of the perturbation differs. 

First, the simplest Method~I of phase reduction is formulated as follows.
For an oscillator described by~$\dot{\vec{x}} = \vec{F}(\vec{x})$ with~$\vec{x} \in {\mathbb R}^d$ with an exponentially stable limit-cycle solution~$\vec{x}_0(t)$ of period~$T$ and frequency~$\omega = 2\pi/T$, a phase~$\vartheta = \omega t$ ($0 \leq t \leq T$,~$0 \leq \vartheta \leq 2\pi$) is introduced on the limit cycle, and the oscillator state on the limit cycle is denoted as~$\vec{\chi}(\vartheta) = \vec{x}_0(\vartheta / \omega)$. 
The basin of attraction of the limit cycle is denoted as~$B \subseteq {\mathbb R}^d$. 
An asymptotic phase function,~$\Theta : B \to [0, 2\pi]$, is introduced such that~$\Theta(\vec{\chi}(\vartheta)) = \vartheta$ on the limit cycle and~$\vec{F}(\vec{x}) \cdot \vec{\nabla} \Theta(x) = \omega$ for all~$\vec{x}$ in~$B$.
The level sets~$I(\vartheta)$ of~$\Theta(\vec{x})$ are called isochrones, which are characterized by a phase~$\vartheta$, and the basin~$B$ is foliated by such isochrones.
Now, for a weakly perturbed oscillator,~$\dot{\vec{x}} = \vec{F}(\vec{x}) + \varepsilon \vec{p}(\vec{x}, t)$, where~$\vec{p}(\vec{x}, t)$~represents the perturbation, the phase~$\vartheta = \Theta(\vec{x})$ of the oscillator obeys~$\dot{\vartheta} = \dot{\Theta}(\vec{x}) = \dot{\vec{x}} \cdot \vec{\nabla} \Theta(\vec{x}) = ( \vec{F}(\vec{x}) + \varepsilon \vec{p}(\vec{x}, t) ) \cdot \vec{\nabla} \Theta(\vec{x}) = \omega + \varepsilon \vec{\nabla} \Theta(\vec{x}) \cdot \vec{p}(\vec{x}, t)$. Though this equation is not closed in~$\theta$, if~$\varepsilon$ is sufficiently small, the oscillator state~$\vec{x}$ is close to a state~$\vec{\chi}(\vartheta)$ on the limit cycle. Then, by approximating~$\vec{x}$ as~$\vec{\chi}(\vartheta)$, we obtain a phase equation~$\dot{\vartheta} = \omega + \varepsilon \vec{Z}(\vartheta) \cdot \vec{p}(\vec{x}_0(\vartheta), t)$ closed in~$\vartheta$ up to~$O(\varepsilon)$, where~$\vec{Z}(\vartheta) = \vec{\nabla} \Theta(\vec{x})\rvert_{\vec{x} = \vec{\chi}(\vartheta)}$ is called the phase sensitivity function (or infinitesimal phase resetting curve).

The above argument is given in the footnote of~\cite{winfree1967biological}, which is considered to be the first derivation of the phase equation for a weakly perturbed oscillator. This method has been extensively used in deriving coupled phase equations from coupled limit-cycle oscillators, including the well-known Kuramoto model exhibiting collective synchronization. 
For example, when there are~$n$ coupled limit-cycle oscillators described as~$\dot{\vec{x}}_k = \vec{F}(\vec{x}_k) + \varepsilon \vec{G}_k(\vec{x}_1, \dotsc, \vec{x}_n)$ for~$k=1, \dotsc, n$, the reduced coupled phase equations are given by~$\dot \theta_k = \omega + \varepsilon \vec{Z}(\theta_k) \cdot \vec{G}_k(\vec{\chi}(\theta_1), \dotsc, \vec{\chi}(\theta_n))$. Typically, by further performing averaging, non-resonant terms are removed, resulting in a phase equation of the form~$\dot \theta_k = \omega + \varepsilon q_k(\theta_1 - \theta_k, \dotsc, \theta_n-\theta_k)$, where the same~$\theta_k$ is used for representing the phase variable after the near-identity transformation. 
We note that, if the coupling~$\vec{G}_k$ is a linear function of~$\vec{x}_1, \dotsc, \vec{x}_n$, then it can be decomposed into terms including only a single variable, e.g.,~$\vec{x}_j$, resulting in a pairwise phase coupling term~$\vec{Z}(\theta_k) \cdot \vec{\chi}(\theta_j)$ or~$q_k(\theta_j-\theta_k)$. Thus, from coupled limit-cycle oscillators with linear pairwise interactions, only pairwise phase models are derived by the Method~I.

Though often overlooked, Kuramoto also developed a more systematic higher-order phase reduction theory in~\cite{Kuramoto_book}, which is called Method~II of phase reduction.
In the Method~I, the geometrically complex asymptotic phase and isochrones were used for the description and deformation of the trajectory due to perturbation was not considered. In Method~II, rather than the curved isochrones, the Floquet coordinates along the limit cycle were employed, and a systematic perturbation theory was formulated to calculate the deformation of the trajectory and the associated phase equation.

Specifically, for a single oscillator of the form~$\dot{\vec{x}} = \vec{F}(\vec{x}) + \varepsilon \vec{G}(\vec{x})$, the trajectory and the associated phase dynamics are assumed to take the form~$\vec{x} = \vec{\chi}(\vartheta) + \varepsilon \vec{\rho}(\vartheta)$ and~$\dot{\vartheta} = \omega + \varepsilon q(\vartheta)$. Here, the phase~$\vartheta$ of the state~$\vec x$ is chosen so that~$\vec{\rho}(\vartheta)$ does not contain the zero Floquet eigenvector along the original limit cycle, i.e.,~$\vec{x}$ and~$\vec{\chi}(\vartheta)$ are on the same tangent space~$T(\vartheta)$ of the isochrone~$I(\vartheta)$. Then, these equations are plugged into the original equation for~$\vec{x}$, yielding an equation that should be satisfied by~$\vec{\rho}$ and~$q(\vartheta)$. Explicitly, it is given in the form 
\begin{align*}
( \omega + \varepsilon q(\vartheta) ) \left[ \frac{d\vec{\chi}(\vartheta)}{d\vartheta} + \varepsilon \frac{d\vec{\rho}(\vartheta)}{d\vartheta} \right]
=
{\vec{F}}(\vec{\chi}(\vartheta) + \varepsilon \vec{\rho}(\vartheta)) + \varepsilon {\vec{G}}(\vec{\chi}(\vartheta) + \varepsilon \vec{\rho}(\vartheta)),
\end{align*}
which is exactly the same as the conjugacy equation~\eqref{eq:ConjEq} by noting the relation~$\vec{e}(\vartheta) = \vec{\chi}(\vartheta) + \vec{\rho}(\vartheta)$.
Then,~$q(\vartheta)$ and~$\vec{\rho}(\vartheta)$ are expanded as~$q = q^{(0)} + \varepsilon q^{(1)} + \cdots$ and~$\vec{\rho} = \vec{\rho}^{(0)} + \varepsilon \vec{\rho}^{(1)} + \dotsb$ and plugged into the above equation, yielding a set of equations.
By construction, those equations are also exactly the same as the homological equations~\eqref{eq:lHom}.
Those equations are further solved by expanding them into the Floquet eigenvectors, yielding explicit expressions for the correction to the orbit and the associated phase equation at each order of~$\varepsilon$.

Thus, the Kuramoto's Method~II and the Parametrization Method in this paper are exactly the same for a single oscillator. In~\cite{Kuramoto_book}, the method is extended to include spatial derivatives to address spatiotemporal dynamics of oscillatory media, resulting in the Kuramoto--Sivashinsky equation obtained at the second order perturbation.

This Method~II has not been used since then to analyze coupled oscillators. Only recently, in~\cite{Kuramoto2019}, Kuramoto extended this method to a discrete network of coupled oscillators.
Similarly to the single-oscillator case, for~$n$ coupled oscillators of the form~$\dot {\vec{x}}_k = \vec{F}(\vec{x}_k) + \varepsilon \vec{G}_k(\vec{x}_1, \dotsc, \vec{x}_n)$ for~$k=1, \dotsc, n$, the oscillator states and their phase dynamics are assumed to be~$\vec{x}_k = \vec{\chi}(\theta_k) + \varepsilon \vec{\rho}(\theta_k)$ and~$\dot{\theta}_k = \omega + \varepsilon q(\theta_k)$, and plugged into the original equations. Then, the conjugacy equations are obtained in the form
\begin{align*}
&
\left( \omega + \varepsilon q_k( \vec \theta) \right) \frac{d{\vec{\chi}}(\theta_k)}{d\theta_k} 
+ \varepsilon \sum_{j=1}^n \left( \omega + \varepsilon q_j(\vec \theta) \right) \frac{\pa {\vec{\rho}}_k( \vec \theta)}{\pa \theta_j} 
\cr
&\qquad=
\vec{F}({\vec{\chi}}(\theta_k) + \varepsilon {\vec{\rho}}_k( \vec \theta)) 
+ \varepsilon \vec{G}_k({\vec{\chi}}(\theta_1) + \varepsilon {\vec{\rho}}_1( \vec \theta),\dotsc, {\vec{\chi}}(\theta_n) + \varepsilon {\vec{\rho}}_n(\vec \theta)),
\end{align*}
which is also the same as the conjugacy equation~\eqref{eq:ConjEq} in this paper. 
When~$\rho_k$ and~$q_k$ are expanded in series of~$\varepsilon$, the equation at each order yields the same homological equation as in this paper.

In Kuramoto's Method II~\cite{Kuramoto2019}, a slightly different method is used in the perturbation expansion.
To simplify the analysis, an additional assumption is introduced, that is, the oscillators are nearly identical and the phase differences are slow variables of~$O(\varepsilon)$. 
Then, the independent variables for the oscillator~$k$ are chosen as the raw phase~$\theta_k$ and the phase differences~$\psi_j = \theta_j - \theta_k$ for~$j=1, \dotsc, n, j \neq k$. 
The equations at~$O(\varepsilon^0)$ and~$O(\varepsilon^1)$ are basically the same as the homological equations in this paper, yielding the same results.
In particular, the $O(\varepsilon^1)$ equation gives the same lowest-order phase equation as the Method~I.
However, the equation at~$O(\varepsilon^2)$ is different from the one in this paper.
In contrast, in the Parametrization Method, no such additional assumption of nearly identical oscillators is made, and the~$O(\varepsilon^2)$ homological equations are directly treated. 

Thus, for some specific class of limit-cycle oscillators, the two method might yield different phase equations at~$O(\varepsilon^2)$.
It would be an interesting future topic to perform actual derivation of the phase equations using both methods and compare whether they yield the same higher-order phase equation or not.
However, it should also be noted that, since the original coupled limit-cycle oscillators are the same, they should eventually describe the same physical phenomena, even if the~$O(\varepsilon^2)$ equations take slightly different forms.

\subsection*{First order phase reduction for pairwise coupling}

Let us consider the following system of three linearly coupled Stuart--Landau oscillators of the form~$\dot {\vec z} = \vec F+\varepsilon\vec G(z)$ given by
\begin{equation}\label{eq:SL_methods}
   \begin{split}
     \dot{z}_1&=(a+ib)z_1+(c+id)\lvert z_1 \rvert^2z_1+\varepsilon \mathrm{e}^{i\varrho}(w_{11}z_1+w_{12}z_2+w_{13}z_3), \\
     \dot{z}_2&=(a+ib)z_2+(c+id)\lvert z_2 \rvert^2z_2+\varepsilon \mathrm{e}^{i\varrho}(w_{21}z_1+w_{22}z_2+w_{23}z_3), \\
     \dot{z}_3&=(a+ib)z_3+(c+id)\lvert z_3 \rvert^2z_3+\varepsilon \mathrm{e}^{i\varrho}(w_{31}z_1+w_{32}z_2+w_{33}z_3),
 \end{split}
\end{equation}
whose parameters are as in the Results, and compute the phase reduction at the first order. 
Using the phase reduction method through parameterization outlined above, we proceed by computing the tangent component (cf.~\eqref{eq:ConjSolve} with the inhomogeneity~\eqref{eq:homological1}), i.e., \begin{equation}\label{eq:tangent_AppA}
\partial_{\vec{\omega}} \ord{\vec{e}}{1}
+
\ord{\vec{f}}{1}
=
\mathbf{T}^{+}
\circ
\Pi
\circ
(\vec{G}\circ\ord{\vec{e}}{0}),
\end{equation}
where 
\begin{align*}
    \mathbf{T} &= \diag(iR\mathrm{e}^{i\theta_j}) & 
    \mathbf{T}^{+}&= \diag((iR)^{-1}\mathrm{e}^{-i\theta_j}) \\
    \Pi &=\diag(\mathrm{e}^{i\theta_j}\Pi(0)\mathrm{e}^{-i\theta_j}) &
    \Pi(0)(x+iy) &=i\left(y-\frac{d}{c}x\right)
\end{align*}
as in the running example of above.

We start by computing 
\begin{align}\label{eq:RHSOrdZero}
    \vec{G}\circ\ord{\vec{e}}{0}=\begin{bmatrix}
        w_{11}R\mathrm{e}^{i(\varrho+\theta_1)}+w_{12}R\mathrm{e}^{i(\varrho+\theta_2)}+w_{13}R\mathrm{e}^{i(\varrho+\theta_3)} \\  w_{21}R\mathrm{e}^{i(\varrho+\theta_1)}+w_{22}R\mathrm{e}^{i(\varrho+\theta_2)}+w_{23}R\mathrm{e}^{i(\varrho+\theta_3)} \\
         w_{31}R\mathrm{e}^{i(\varrho+\theta_1)}+w_{32}R\mathrm{e}^{i(\varrho+\theta_2)}+w_{33}R\mathrm{e}^{i(\varrho+\theta_3)}
    \end{bmatrix}.
\end{align}
We then apply the projection~$\Pi$ and then~$\mathbf{T}^{+}$ to obtain
\begin{equation}\label{eq:f1_expl_AppA}
\vec{f}^{(1)} =
\begin{bmatrix}
w_{11}\sin(\varrho)
+ w_{12}\sin(\varrho+\theta_2-\theta_1)
+ w_{13}\sin(\varrho+\theta_3-\theta_1) \\
\quad - \frac{d}{c} \Big(
w_{11}\cos(\varrho)
+ w_{12}\cos(\varrho+\theta_2-\theta_1)
+ w_{13}\cos(\varrho+\theta_3-\theta_1)
\Big) \\[0.8em]

w_{21}\sin(\varrho+\theta_1-\theta_2)
+ w_{22}\sin(\varrho)
+ w_{23}\sin(\varrho+\theta_3-\theta_2) \\
\quad - \frac{d}{c} \Big(
w_{21}\cos(\varrho+\theta_1-\theta_2)
+ w_{22}\cos(\varrho)
+ w_{23}\cos(\varrho+\theta_3-\theta_2)
\Big) \\[0.8em]

w_{31}\sin(\varrho+\theta_1-\theta_3)
+ w_{32}\sin(\varrho+\theta_2-\theta_3)
+ w_{33}\sin(\varrho) \\
\quad - \frac{d}{c} \Big(
w_{31}\cos(\varrho+\theta_1-\theta_3)
+ w_{32}\cos(\varrho+\theta_2-\theta_3)
+ w_{33}\cos(\varrho)
\Big)
\end{bmatrix}.
\end{equation}
By using the identity~$A\cos(t)+B\sin(t)=\sqrt{A^2+B^2}\sin(t+\arctan(\frac{A}{B}))$, we obtain
\begin{equation}\label{eq:f1_comp_AppA}
    \vec{f}^{(1)}=K\begin{bmatrix}
        w_{11}\sin(C)+w_{12}\sin(\theta_2-\theta_1+C)+w_{13}\sin(\theta_3-\theta_1+C) \\ w_{21}\sin(\theta_1-\theta_2+C)+w_{22}\sin(C)+w_{23}\sin(\theta_3-\theta_2+C) \\ w_{31}\sin(\theta_1-\theta_3+C)+w_{32}\sin(\theta_2-\theta_3+C)+w_{33}\sin(C)
    \end{bmatrix},
\end{equation}
where~$K=\sqrt{1+\frac{d^2}{c^2}}$ and~$C=\varrho+\arctan(-\frac{d}{c})$. Note that this is already in normal form. 

\subsection*{Second order phase reduction for pairwise coupling: Emergent nonpairwise (EN) interactions}

We now compute the second order phase reduction for linearly coupled SL oscillators~\eqref{eq:SL}. 
To do this, we need to compute the embedding at the first order
\begin{equation}
    \ord{\vec{e}}{1} = \mathbf{T}\ord{\vec{g}}{1}+\mathbf{N}\ord{\vec{h}}{1},
\end{equation} 
where $\ord{\vec{g}}{1}, \ord{\vec{h}}{1}$ are determined by~\eqref{eq:ConjSolve}. 
Applying the corresponding projections~\eqref{eq:SLProjections} and pseudo-inverses~\eqref{eq:SLInverses} to~\eqref{eq:RHSOrdZero}, one obtains that the right hand side~\eqref{eq:ConjSolve} only contains resonant terms.
Thus, we can choose $\ord{\vec{g}}{1}=0$ (as~$\vec{f}^{(1)}$ is already in normal form as noted above) and we have $\partial_{\vec{\omega}}\ord{\vec{h}}{1}=0$ (cf.~\cite[Eq.~(39)]{Bick2025}).
Inverting $L=\diag(-2a)$, we obtain
\begin{equation}\label{eq:h1_AppB}
    \ord{\vec{h}}{1}=\frac{R}{2ac}\begin{bmatrix}
        \tilde{h}_{11} \\ \tilde{h}_{12} \\ \tilde{h}_{13}
    \end{bmatrix},
\end{equation} 
where 
\begin{equation}\label{eq:htilde_AppB}
    \begin{split}
         \tilde{h}_{11} =w_{11}\cos(\varrho)+w_{12}\cos(\varrho+\theta_2-\theta_1)+w_{13}\cos(\varrho+\theta_3-\theta_1), \\ \tilde{h}_{12}=w_{21}\cos(\varrho+\theta_1-\theta_2)+w_{22}\cos(\varrho)+w_{23}\cos(\varrho+\theta_3-\theta_2),  \\ \tilde{h}_{13}=w_{31}\cos(\varrho+\theta_1-\theta_3)+w_{32}\cos(\varrho+\theta_2-\theta_3)+w_{33}\cos(\varrho). 
    \end{split}
\end{equation}
Hence, we we have the first-order term for the embedding
\begin{equation}\label{eq:e1_AppB}
    \ord{\vec{e}}{1}=\frac{(c+id)R}{2ac}\begin{bmatrix}
        \mathrm{e}^{i\theta_1}\tilde{h}_{11} \\  \mathrm{e}^{i\theta_2}\tilde{h}_{12} \\  \mathrm{e}^{i\theta_3}\tilde{h}_{13}
    \end{bmatrix}.
\end{equation}

With the first-order embedding at hand, we can now compute the second order phase interaction from the inhomogeneity~\eqref{eq:homological2}.
As in~\eqref{eq:InhomMain} in the main text, we denote the three terms entering the second order inhomogeneity by
\begin{subequations}
\begin{align}
\label{eq:InhomU}
\vec U &= \vec G'(\ord{\vec e}{0}) \cdot  \ord{\vec e}{1}, \\
\label{eq:InhomV}
\vec V &= -(\ord{\vec e}{1})' \cdot \ord{\vec f}{1}, \\
\label{eq:InhomW}
\vec W &= \frac{1}{2}\vec F''(\ord{\vec e}{0})\left[\ord{\vec e}{1}, \ord{\vec e}{1}\right].
\end{align}
\end{subequations}
so that
\[
\vec\eta^{(2)}
=
\vec U+\vec V +\vec W.
\]
In the following, we compute the projections of the three terms
separately.

\paragraph{First term~$\vec{U}$}
For the first term~$\vec{U}$ of the inhomogeneity~\eqref{eq:InhomU}, we have
\begin{displaymath}
    \vec U = \frac{(c+id)}{2ac}R\begin{bmatrix}w_{11}\mathrm{e}^{i\theta_1}\tilde{h}_{11}+w_{12}\mathrm{e}^{i\theta_2}\tilde{h}_{12}+w_{13}\mathrm{e}^{i\theta_3}\tilde{h}_{13} \\ w_{21}\mathrm{e}^{i\theta_1}\tilde{h}_{11}+w_{22}\mathrm{e}^{i\theta_2}\tilde{h}_{12}+w_{23}\mathrm{e}^{i\theta_3}\tilde{h}_{13} \\ w_{31}\mathrm{e}^{i\theta_1}\tilde{h}_{11}+w_{32}\mathrm{e}^{i\theta_2}\tilde{h}_{12}+w_{33}\mathrm{e}^{i\theta_3}\tilde{h}_{13}    
    \end{bmatrix}
\end{displaymath}
Thus, we have
\begin{equation}\label{eq:finalT1_AppB}
    \mathbf{T}^{+}\circ \Pi \circ \vec U = \frac{1}{2a}\Big(1+\frac{d^2}{c^2}\Big)\begin{bmatrix}
        w_{12}\tilde{h}_{12}\sin(\theta_2-\theta_1+\varrho)+w_{13}\tilde{h}_{13}\sin(\theta_3-\theta_1+\varrho) \\ w_{21}\tilde{h}_{11}\sin(\theta_1-\theta_2+\varrho)+w_{23}\tilde{h}_{13}\sin(\theta_3-\theta_2+\varrho) \\ w_{31}\tilde{h}_{11}\sin(\theta_1-\theta_3+\varrho)+w_{32}\tilde{h}_{12}\sin(\theta_2-\theta_3+\varrho)
    \end{bmatrix},
\end{equation}
with~$\tilde{h}_{1j}$ as defined in Eq.~\eqref{eq:htilde_AppB}.

\allowdisplaybreaks

\paragraph{Second term~$\vec{V}$}
We first compute the second term~$\vec{V}$ of the inhomogeneity~\eqref{eq:InhomV}.
First, the Jacobian of the first order of the embedding~$\ord{\vec{e}}{1}$ is
\begin{displaymath}
(\ord{\vec{e}}{1})' = \frac{(c+id)R}{2ac}
\begin{bmatrix}
\mathrm{e}^{i\theta_1}(i\tilde{h}_{11}+\alpha_{11}) & \mathrm{e}^{i\theta_1}\alpha_{12} & \mathrm{e}^{i\theta_1}\alpha_{13} \\
\mathrm{e}^{i\theta_2}\alpha_{21} & \mathrm{e}^{i\theta_2}(i\tilde{h}_{12}
+ \alpha_{22}) & \mathrm{e}^{i\theta_2}\alpha_{23} \\
\mathrm{e}^{i\theta_3}\alpha_{31} & \mathrm{e}^{i\theta_3}\alpha_{32} & \mathrm{e}^{i\theta_3}(i\tilde{h}_{13}
+ \alpha_{33})
\end{bmatrix},
\end{displaymath}
where
\begin{displaymath}
\small
\begin{aligned}
\alpha_{11} &= w_{12}\sin(\varrho+\theta_2-\theta_1)
+ w_{13}\sin(\varrho+\theta_3-\theta_1),\\
\alpha_{12} &= - w_{12} \sin(\varrho+\theta_2-\theta_1),\\
\alpha_{13} &= - w_{13}  \sin(\varrho+\theta_3-\theta_1),\\
\alpha_{21} &= - w_{21} \sin(\varrho+\theta_1-\theta_2),\\
\alpha_{22} &= w_{21}\sin(\varrho+\theta_1-\theta_2)
+ w_{23}\sin(\varrho+\theta_3-\theta_2),\\
\alpha_{23} &= - w_{23}  \sin(\varrho+\theta_3-\theta_2),\\
\alpha_{31} &= - w_{31}  \sin(\varrho+\theta_1-\theta_3),\\
\alpha_{32} &= - w_{32} \sin(\varrho+\theta_2-\theta_3),\\
\alpha_{33} &= w_{31}\sin(\varrho+\theta_1-\theta_3)
+ w_{32}\sin(\varrho+\theta_2-\theta_3).
\end{aligned}
\normalsize
\end{displaymath}
With~$\vec{f}^{(1)}$ and~$K$ as in~\eqref{eq:f1_comp_AppA} above,
we have
\begin{displaymath}
    \vec V = -\frac{(c+id)KR}{2ac}\begin{bmatrix}
        \mathrm{e}^{i\theta_1}\Big((i\tilde{h}_{11}
+\alpha_{11})f_1^{(1)}+\alpha_{12}f_2^{(1)}+\alpha_{13}f_3^{(1)}\Big) \\ \mathrm{e}^{i\theta_2}\Big(\alpha_{21}f_1^{(1)}+(i\tilde{h}_{12}
+ \alpha_{22})f_2^{(1)}+\alpha_{23}f_3^{(1)}\Big) \\\mathrm{e}^{i\theta_3}\Big(\alpha_{31}f_1^{(1)}+\alpha_{32}f_2^{(1)}+(i\tilde{h}_{13}+\alpha_{33})f_3^{(1)}\Big)
    \end{bmatrix}.
\end{displaymath}
With the shorthand
\begin{displaymath}
    \begin{split}        \mathcal{Y}_1=\alpha_{11}f_1^{(1)}+\alpha_{12}f_2^{(1)}+\alpha_{13}f_3^{(1)}, \\   \mathcal{Y}_2=\alpha_{21}f_1^{(1)}+\alpha_{22}f_2^{(1)}+\alpha_{23}f_3^{(1)}, \\   \mathcal{Y}_3=\alpha_{31}f_1^{(1)}+\alpha_{32}f_2^{(1)}+\alpha_{33}f_3^{(1)},
    \end{split}
\end{displaymath} 
we obtain
\begin{displaymath}
    \vec V = - \frac{(c+id)KR}{2ac}\begin{bmatrix}
        \mathrm{e}^{i\theta_1}\Big(i\tilde{h}_{11}
f_1^{(1)}+\mathcal{Y}_1 \Big) \\ \mathrm{e}^{i\theta_2}\Big(i\tilde{h}_{12}
f_2^{(1)}+\mathcal{Y}_2\Big) \\\mathrm{e}^{i\theta_3}\Big(i\tilde{h}_{13}f_3^{(1)}+\mathcal{Y}_3\Big)
    \end{bmatrix}.
\end{displaymath}

This allows us to compute the resulting second order contributions.
By applying the projection, we find that
\begin{equation}\label{eq:finalT3_AppB}
    \mathbf{T}^{+}\circ \Pi \circ \vec V =-\frac{K}{2ac}\begin{bmatrix}
        (c+\frac{d^2}{c})\tilde{h}_{11}
f_1^{(1)}  \\  (c+\frac{d^2}{c})\tilde{h}_{12}
f_2^{(1)}  \\  (c+\frac{d^2}{c})\tilde{h}_{13}
f_3^{(1)}
    \end{bmatrix}.
\end{equation}

\paragraph{Third term~$\vec{W}$}
Finally, we compute the third term~$\vec{W}$ of the inhomogeneity~\eqref{eq:InhomW} for general parameters of the SL system.
Recall that
\begin{align*}
\vec{F} &=
\begin{bmatrix}
(a+ib)z_1+(c+id)\lvert z_1 \rvert^2z_1 \\
(a+ib)z_2+(c+id)\lvert z_2 \rvert^2z_2 \\
(a+ib)z_3+(c+id)\lvert z_3 \rvert^2z_3
\end{bmatrix},
\end{align*}
as well as
\begin{align*}
\ord{\vec{e}}{0} &=
\begin{bmatrix}
R \mathrm{e}^{i\theta_1} \\
R \mathrm{e}^{i\theta_2} \\
R \mathrm{e}^{i\theta_3}
\end{bmatrix},
&
\ord{\vec{e}}{1}
= \frac{(c+id)R}{2ac}
\begin{bmatrix}
\mathrm{e}^{i\theta_1}\tilde{h}_{11} \\
\mathrm{e}^{i\theta_2}\tilde{h}_{12} \\
\mathrm{e}^{i\theta_3}\tilde{h}_{13}
\end{bmatrix}.
\end{align*}

We compute the second derivative of~$\vec{F}$.
The Jacobian of one component of~$\vec{F}$ evaluates to 
\begin{displaymath}
    F'(z) =\begin{bmatrix} \frac{\partial F}{\partial z} \\\frac{\partial F}{\partial \bar{z}}\end{bmatrix}= \begin{bmatrix} (a+ib)+2(c+id)z\bar{z} \\ (c+id)z^2\end{bmatrix}.
\end{displaymath}
Taking derivatives again to obtain the Hessian, we have
\begin{displaymath}
   F''(z) = 
\begin{bmatrix}
\frac{\partial^2 F}{\partial z^2} & \frac{\partial^2 F}{\partial \bar{z} \partial z} \\
\frac{\partial^2 F}{\partial z \partial \bar{z}} & \frac{\partial^2 F}{\partial \bar{z}^2}
\end{bmatrix}
=
\begin{bmatrix}
2(c+id)\bar{z} & 2(c+id)z \\
2(c+id)z & 0
\end{bmatrix}.
\end{displaymath}

Given a complex number~$v$, the bilinear form~$F''(z)[v, v]$ evaluates to
\begin{displaymath}
 F''(z)[v, v]=  \Bigg( \begin{bmatrix} v & \bar{v} \end{bmatrix}\cdot
\begin{bmatrix} 2(c+id) \bar{z} & 2(c+id)z \\ 2(c+id)z & 0 \end{bmatrix}\Bigg)\cdot
\begin{bmatrix} v \\ \bar{v} \end{bmatrix}=2(c+id)\bar{z} v^2 + 4(c+id)z \abs{v}^2.
\end{displaymath}
For one component of the embedding~$\ord{\vec{e}}{0}_j$, the Hessian is
\begin{displaymath}
   F''(\ord{\vec{e}}{0}_j) = R 
\begin{bmatrix}
2(c+id)\mathrm{e}^{-i\theta_j} & 2(c+id)\mathrm{e}^{i\theta_j} \\
2(c+id)\mathrm{e}^{i\theta_j} & 0
\end{bmatrix},
\end{displaymath}
so that the~$j$th component of the bilinear form is given by
\[
\left(
\frac{1}{2}
F''(\ord{\vec{e}}{0})
\left[
\ord{\vec{e}}{1},
\ord{\vec{e}}{1}
\right]
\right)_j
=
(c+id)R \mathrm{e}^{-i\theta_j}
\left(\ord{\vec{e}}{1}_j\right)^2
+
2(c+id)R \mathrm{e}^{i\theta_j}
\abs{\ord{\vec{e}}{1}_j}.
\]
Treating the square and square modulus of the first order embedding separately, we have 
\begin{displaymath}
    \Big(\ord{\vec{e}}{1}_j\Big)^2
    =\Bigg(\frac{(c+id)R}{2ac}
\mathrm{e}^{i\theta_j}\tilde{h}_{1j}\Bigg)^2=\frac{(c+id)^2R^2}{4a^2c^2}
\mathrm{e}^{2i\theta_j}\left(\tilde{h}_{1j}\right)^2,
\end{displaymath}
and 
\[
\big\lvert\ord{\vec{e}}{1}_j \big\rvert^2
=
\left\lvert
\frac{(c+id)R}{2ac}
\mathrm{e}^{i\theta_j}\tilde{h}_{1j}
\right\rvert^2
=
\frac{(c^2+d^2)R^2}{4a^2c^2}
\lvert \tilde{h}_{1j} \rvert^2
\]
So that 
\begin{displaymath}
   \Bigg( \frac{1}{2} F''(\ord{\vec{e}}{0})\big[ \ord{\vec{e}}{1}, \ord{\vec{e}}{1}\big]\Bigg)_j=[2(c^2+d^2)+(c+id)^2]\frac{R^3}{4a^2c^2}\mathrm{e}^ {i\theta_j}
\left(\tilde{h}_{1j}\right)^2.
\end{displaymath}

Written as a vector of three components, we have 
\begin{displaymath}
    \vec W = \frac{1}{2} F_0{}''(\ord{\vec{e}}{0})\big[ \ord{\vec{e}}{1}, \ord{\vec{e}}{1}\big]=\frac{(3c^2+d^2+2icd)R^3}{2a^2c^2}\begin{bmatrix}
        \mathrm{e}^{i\theta_1}\left(\tilde{h}_{11}\right)^2 \\
         \mathrm{e}^{i\theta_2}\left(\tilde{h}_{12}\right)^2 \\
        \mathrm{e}^{i\theta_3}\left(\tilde{h}_{13}\right)^2
    \end{bmatrix}.
\end{displaymath}
so that with the projections we obtain
\begin{equation}\label{eq:finalT2_AppB}
    \mathbf{T}^{+}\circ \Pi \circ \vec W =-d\left(5c-\frac{d^2}{c}\right)\left(\frac{R}{ac}\right)^2\begin{bmatrix}
        \left(\tilde{h}_{11}\right)^2 \\
         \left(\tilde{h}_{12}\right)^2 \\
        \left(\tilde{h}_{13}\right)^2
    \end{bmatrix}.
\end{equation} 

As we will not proceed to compute the third order, we do not need to compute the second order of the embedding. Moreover, in this work we will not compute the second order for PN interactions, i.e., hypergraphs, which is left for future works.

\subsection*{Explicit expression for the first component of the
second order linear phase reduction}

The second order phase reduction is determined by 
$\mathbf{T}^+\circ\Pi\circ
\big(
\vec U+\vec V +\vec W
\big)$.
More precisely, the first component~$f^{(2)}
_1=\big[\vec 
f^{(2)}
\big]
_1
$ of the second order contribution to the phase dynamics is
\begin{equation}
f^{(2)}_1
=
\left[
\mathbf{T}^+\circ\Pi\circ\vec U
\right]_1
+
\left[
\mathbf{T}^+\circ\Pi\circ\vec V
\right]_1
+
\left[
\mathbf{T}^+\circ\Pi\circ\vec W
\right]_1.
\label{eq:explicit_second_order_decomposition}
\end{equation}
We now collect the terms computed above and write them out explicitly.
We focus on the case of straight
isochrones, i.e.,~$d=0$, so that~$K=1$ and~$C=\varrho$, and assume no
self-coupling, i.e.,
\[
w_{11}=w_{22}=w_{33}=0.
\]

First, for~$\vec U$, we have
\begin{align}
\left[
\mathbf{T}^+\circ\Pi\circ\vec U
\right]_1
=
\frac{1}{2a}
\Big[
&w_{12}\tilde h^{(1)}_2
\sin(\theta_2-\theta_1+\varrho)
\nonumber\\
&+
w_{13}\tilde h^{(1)}_3
\sin(\theta_3-\theta_1+\varrho)
\Big].
\end{align}
Developing the above expression, we obtain
\begin{align}
\left[
\mathbf{T}^+\circ\Pi\circ\vec U
\right]_1
=
\frac{1}{4a}
\Big[
&w_{12}w_{21}C_{1221}
+w_{12}w_{23}C_{1223}
\nonumber\\
&+
w_{13}w_{31}C_{1331}
+w_{13}w_{32}C_{1332}
\Big],
\label{eq:explicit_U_first_component}
\end{align}
where
\begin{subequations}
\begin{align}
C_{1221}
&=
\sin(2\varrho)
+
\sin\left(2(\theta_2-\theta_1)\right),
\\
C_{1223}
&=
\sin(\theta_3-\theta_1+2\varrho)
+
\sin(2\theta_2-\theta_3-\theta_1),
\\
C_{1331}
&=
\sin(2\varrho)
+
\sin\left(2(\theta_3-\theta_1)\right),
\\
C_{1332}
&=
\sin(\theta_2-\theta_1+2\varrho)
+
\sin(2\theta_3-\theta_2-\theta_1).
\end{align}
\end{subequations}


Second, for~$\vec V$ we have
\begin{equation}
\left[
\mathbf{T}^+\circ\Pi\circ\vec V
\right]_1
=
-\frac{1}{2a}
\tilde h^{(1)}_1 f^{(1)}_1.
\end{equation}
Developing the above expression, we obtain
\begin{equation}
\left[
\mathbf{T}^+\circ\Pi\circ\vec V
\right]_1
=
-\frac{1}{4a}
\left[
w_{12}^2D_{1212}
+2w_{12}w_{13}D_{1213}
+w_{13}^2D_{1313}
\right],
\label{eq:explicit_V_first_component}
\end{equation}
where
\begin{subequations}
\begin{align}
D_{1212}
&=
\sin\left(2(\theta_2-\theta_1)+2\varrho\right),
\\
D_{1213}
&=
\sin(\theta_2+\theta_3-2\theta_1+2\varrho),
\\
D_{1313}
&=
\sin\left(2(\theta_3-\theta_1)+2\varrho\right).
\end{align}
\end{subequations}


Third, setting~$d=0$ we find that the tangential projection of~$\vec W$ vanishes
\begin{equation}
\left[
\mathbf{T}^+\circ\Pi\circ\vec W
\right]_1
=
0.
\label{eq:explicit_W_first_component}
\end{equation}

Combining
Eqs.~\eqref{eq:explicit_U_first_component},
\eqref{eq:explicit_V_first_component}, and
\eqref{eq:explicit_W_first_component}, we finally obtain the expression explicitly presented in the Results, namely,
\begin{align}
f^{(2)}_1
=
\frac{1}{4a}
\Big[
&w_{12}w_{21}C_{1221}
+w_{12}w_{23}C_{1223}
+w_{13}w_{31}C_{1331}
+w_{13}w_{32}C_{1332}
\nonumber\\
&-
w_{12}^2D_{1212}
-2w_{12}w_{13}D_{1213}
-w_{13}^2D_{1313}
\Big].
\label{eq:explicit_second_order_first_component}
\end{align}
Thus, the phase dynamics of oscillator~1 (and with the appropriate adjustment of indices analogously also for the other component) up to second order in~$\varepsilon$ is
\begin{equation}
\dot{\theta}_1=
\omega
+\varepsilon
\left[
w_{12}\sin(\varrho+\theta_2-\theta_1)
+w_{13}\sin(\varrho+\theta_3-\theta_1)
\right]+
\varepsilon^2
f^{(2)}_1
+O(\varepsilon^3).
\label{eq:explicit_phase_dynamics_first_component}
\end{equation}

\subsection*{First order phase reduction for nonpairwise coupling: Physical nonpairwise (PN) interactions}
For what concerns nonpairwise coupling, i.e., physical nonpairwise (PN) interactions, not every possible coupling gives a phase model in normal form. 

We now discuss \emph{all} resonant cubic coupling terms for the network of three Stuart--Landau oscillators~\eqref{eq:SL}, and to clarify why only some of them yield physical nonpairwise (PN) interactions, which are the ones used in the Results. 

The phase reduced model always takes a~$-\theta_j$ in the~$\sin$ of the coupling function relative to the~$j$th~oscillator. 
Hence, from a general coupling function of the form~$z_1^{n_1}\dotsb z_{j-1}^{n_{j-1}}z_j^{n_j}z_{j+1}^{n_{j+1}}\dotsb z_l^{n_l}$, the coupling in the phase reduced model will be of the form \begin{displaymath}
    \dot{\theta}_j\propto \sin(n_1\theta_1+n_{j-1}\theta_{j-1}+n_j\theta_j+n_{j+1}\theta_{j+1}+\dotsb+n_l\theta_l).
\end{displaymath}
The condition to have resonant terms is then~$n_1+\dotsb+n_{j-1}+n_j+n_{j+1}+\dotsb+n_l=0$. 
If we call~$n_*$ total order of the conjugated exponents and~$n$~the total order of the non-conjugated ones, we have the relation~$n-n_*=0$. This is a further confirmation that even couplings do not appear in the first order phase model, as shown by León et at.~\cite{leon2025theory}.

Here, we focus on~$3$-body interactions and, more specifically, for the~$3$ oscillators setting we have dealt throughout the whole paper. Because of the resonance condition, the minimal monomial coupling yielding~$3$-body interactions is cubic, which is the one we discuss below. The same reasoning, \textit{mutatis mutandis}, applies to any kind of many-body interaction. \\

Every resonant cubic coupling function for oscillator~$j\in \{1,2,3\}$ has the general form
\begin{equation}
G_j = G_j^{\mathrm{cub}} := \mathrm{e}^{i\xi} \, w_{j p q r}\, z_p z_q \bar z_r,
\qquad p,q,r \in \{1,2,3\},
\label{eq:cubic-general}
\end{equation}
where~$w_{jpqr} = w_{jqpr}$ is symmetric under exchange of the two unconjugated indices, and summation over repeated indices is implied. Using
the same argument as above, the resulting contribution to the first order phase equation for~$\theta_j$ is
\begin{equation}
\dot\theta_j \propto \, w_{jpqr}\,
\sin\!\big(\theta_p + \theta_q - \theta_r - \theta_j + \xi\big).
\label{eq:cubic-phase-general}
\end{equation}
For our network of~$n=3$ oscillators,~$p,q,r$ necessarily coincide with~$j$ or with one of the (at most two) other oscillators, so that Eq.~\eqref{eq:cubic-phase-general} is constrained to a small number of cases, which we enumerate below (note that indices~$k,l$ always denote the two oscillators other than~$j$, i.e.,~$k\neq l$).

\paragraph{Self-coupling}
Let~$p=q=r=j$,~$G_j \propto z_j^2\bar z_j = \lvert z_j\rvert^2 z_j$. On the
limit cycle, this reduces to a constant multiple of~$z_j$ itself, so this
term merely sums to the intrinsic cubic nonlinearity
$(c+id)\lvert z_j\rvert^2 z_j$ already present in~$F_j$ (rotated by
$\rho$). Hence, it carries no information about the network and is not an actual coupling.

\paragraph{Nonlinear pairwise coupling}
For~$\{p,q,r\}$ involving only~$j,k$, two cases arise:
\begin{align}
G_j &\propto z_j^2 \bar z_k \;\;\Rightarrow\;\;
\dot\theta_j \propto  w_{jjjk} \sin(\theta_j-\theta_k+\xi), \label{eq:pairwise_nlna}\\
G_j &\propto z_j z_k \bar z_j = \lvert z_j\rvert^2 z_k \;\;\Rightarrow\;\;
\dot\theta_j \propto  w_{jjkj} \sin(\theta_k-\theta_j+\xi). \label{eq:pairwise_nlnb}
\end{align}
Both produce the same Kuramoto--Sakaguchi harmonic as the linear coupling (Eq.~\eqref{eq:pairwise_nlna} repulsive, Eq.~\eqref{eq:pairwise_nlnb} attractive), weighted by a factor (absorbed into the weights~$w_{jpqr}$) that is constant
on the unperturbed torus. These are pairwise couplings, distinguishable from the linear ones only at second order, or if the oscillators are not identical, which are cases not considered here.

\paragraph{Coupling with no phase dependence}
For~$G_j \propto z_j z_k \bar z_k = \lvert z_k\rvert^2 z_j$, we have
$p+q-r = j+k-k = j$, so Eq.~\eqref{eq:cubic-phase-general} gives
$\dot\theta_j \propto w \sin(\rho)$. Since~$\lvert z_k\rvert^2 = R^2$ is
constant on the unperturbed torus, this term does not depend on any phase
difference at all: At first order, its only effect is a constant shift of the frequency of
oscillator~$j$.

\paragraph{Coupling with no self-dependence}
For~$G_j \propto z_j z_k \bar z_l$, we obtain
$\dot\theta_j \propto w \sin(\theta_k-\theta_l+\xi)$, a term that depends on the relative phase of the other two oscillators but not on~$\theta_j$. Although it involves~$z_j, z_k, z_l$ jointly through the coupling function,
it does not couple~$\theta_j$ to~$\theta_k$ and~$\theta_l$, but to their difference. We do not consider this case further here, though it may be relevant for driven or forced generalizations of the present framework.

\paragraph{Physical nonpairwise (PN) couplings}
The remaining two cases are the couplings used in the Results,
distinguished by which index is conjugated:
\begin{align}
G_j &\propto z_k^2 \bar z_l \;\;\Rightarrow\;\;
\dot\theta_j \propto  w \sin(2\theta_k-\theta_l-\theta_j+\xi),
&& \text{$(2,-1,-1)$, ``asymmetric''}, \label{eq:PN_211}\\
G_j &\propto z_k z_l \bar z_j \;\;\Rightarrow\;\;
\dot\theta_j \propto w \sin(\theta_k+\theta_l-2\theta_j+\xi),
&& \text{$(1,1,-2)$, ``symmetric''}. \label{eq:PN_112}
\end{align}
These are precisely Eqs.~\eqref{eq:reduced_cubic1} and~\eqref{eq:reduced_cubic2} in the main
text, reproduced here as Eqs.~\eqref{eq:PN_211}--\eqref{eq:PN_112}. Unlike
the cases above, both terms couple~$\theta_j$ to the joint configuration of
the other two oscillators, and generate, already at first order, the same
$(2,-1,-1)$ and~$(1,1,-2)$ harmonics that arise as EN at second order from
linear pairwise coupling (cf.\ Eqs.~\eqref{eq:lin_second_C} and
\eqref{eq:lin_second_D}). This is why they are the couplings used for the
coupling design procedure of the Results.

\subsection*{Details on the numerical methods}

In this Section, we briefly describe the numerical methods used for the mixed order phase reduction and the coupling design.

The boundaries shown in Figs.~\ref{fig:mixed_rho_eps}, \ref{fig:mixed_rho_eta}, and~\ref{fig:phase_diagram} are obtained from the linear stability of the two phase-locked states that organize the dynamics of three identical, all-to-all coupled oscillators, namely, the synchronized state~$\theta_1=\theta_2=\theta_3$ and the twisted state~$\theta_j=\theta_1+2\pi(j-1)/3$.
Both states are invariant for all the models considered, i.e., the full Stuart--Landau system~\eqref{eq:SL} with linear and nonlinear coupling, the system subject to coupling design~\eqref{eq:engineered_SL_1}, and the first-, second-, and mixed-order phase reductions~\eqref{eq:1.5reduction}, since they are fixed by the symmetry of the system (permutations of the oscillators combined with phase shifts).
Their loss or gain of stability, therefore, determines the boundaries between the synchronized, bistable, and incoherent regions: The boundary between the synchronized and the bistable regions (dashed lines) corresponds to the value of~$\varrho$ at which the twisted state becomes stable, while the boundary between the bistable and the incoherent regions (solid lines) corresponds to the value of~$\varrho$ at which the synchronized state loses stability.

For the phase models, we compute the Jacobian of the phase dynamics at the two states by finite differences and remove the neutral direction associated with a uniform phase shift.
For the full Stuart--Landau systems, the synchronized and twisted states correspond to rotating waves~$z_j=r\,\mathrm{e}^{i(\Omega t+\theta_j)}$, where the radius~$r$ and the frequency~$\Omega$ are fixed by the coupling.
Since all the couplings considered are invariant under a global rotation~$z_j\mapsto\mathrm{e}^{i\alpha}z_j$, the rotating waves are equilibria in a frame rotating with frequency~$\Omega$, and their stability is determined by the eigenvalues of the Jacobian in such a frame, once the neutral direction associated with the rotation is removed.
For each value of~$\varepsilon$ (Figs.~\ref{fig:mixed_rho_eps} and~\ref{fig:phase_diagram}) or~$\eta$ (Fig.~\ref{fig:mixed_rho_eta}), we evaluate the largest real part of the remaining eigenvalues on a grid in~$\varrho\in[0,\pi]$ and locate its zero crossings by root finding.

Note that the phase diagrams are symmetric under~$\varrho\mapsto-\varrho$, so that it is enough to consider~$\varrho\in[0,\pi]$.
As a consequence, the region in which only the synchronized state is stable and the incoherent region are centred at~$\varrho=0$ and~$\varrho=\pi$, respectively.
For values of~$\varepsilon$ or~$\eta$ beyond which one of these regions shrinks to zero width, the stability of the corresponding state does not change for any~$\varrho$, e.g., once the synchronized region disappears, the twisted state is stable for all~$\varrho$.
The boundaries, therefore, end at~$\varrho=0$ or~$\varrho=\pi$, which we compute by root finding the largest real part of the eigenvalues of the twisted state at~$\varrho=0$ and of the synchronized state at~$\varrho=\pi$, as a function of~$\varepsilon$ or~$\eta$.

Note that a direct classification by simulating the dynamics from a twisted initial condition is not reliable.
In fact, since the twisted state is invariant, a trajectory starting exactly on it leaves only through the amplification of round-off errors, which strongly delays the detection of its instability and shifts the boundary between the synchronized and the bistable regions towards smaller values of~$\varrho$.
We have, nevertheless, checked the results against direct simulations, performed with a Runge--Kutta (4,5) scheme (\texttt{ode45} in MATLAB, relative and absolute tolerances~$10^{-8}$ and~$10^{-10}$), starting from twisted states perturbed by an asymmetric perturbation of amplitude~$10^{-2}$ and integrating up to~$t=400/\varepsilon$.

\section*{Data availability} 

No datasets were generated or analyzed during the current study. All findings are based on numerical simulations described in the Methods section, which can be reproduced using the details provided therein.

\section*{Code availability} 

The code for the numerical simulations presented in this article is available upon reasonable request to the corresponding author.

\section*{Acknowledgments} 

The authors are grateful to Bob Rink, Babette de Wolff, and Iván León for discussions.
R.M. acknowledges JSPS KAKENHI 24KF0211 for financial support. H.N. acknowledges JSPS KAKENHI 25H01468, 25K03081, and 22H00516 for financial support. C.B. acknowledges support by the project ‘BeyondTheEdge: Higher-Order Networks and Dynamics’ (European Union, REA Grant Agreement No. 101120085). R.M. is grateful to VU Amsterdam for the hospitality. C.B. is grateful to the Institute of Science Tokyo and RIKEN iTHEMS for the hospitality.

\section*{Author Contributions}
 
C.B. conceptualized the study. All authors developed the methodology and carried out the theoretical analysis. R.M. carried out the simulations and curated the visualization. R.M. and C.B. wrote the manuscript. All authors reviewed and edited the manuscript.

\section*{Competing Interests} 

The authors declare no competing interests.

\end{document}